\documentclass[12pt]{article}
\usepackage{amsfonts}
\usepackage{bbm}
\usepackage{mathrsfs}
\usepackage{amssymb}

\def\hang{\hangindent\parindent}
\def\textindent#1{\indent\llap{#1\enspace}\ignorespaces}
\def\re{\par\hang\textindent}

\title{Computation of  Minimal Filtered Free
Resolutions over\\ $\mathbb{N}$-Filtered  Solvable Polynomial
Algebras}

\author{Huishi Li\\
{\small Department of Applied Mathematics},
{\small College of Information Science and Technology}\\
{\small Hainan University}, {\small  Haikou 570228, China}}

\date{}

\begin{document}
\maketitle
\begin{center}
\begin{minipage}{140mm}
{\small {\bf Abstract.} Let $A=K[a_1,\ldots ,a_n]$ be a weighted
$\mathbb{N}$-filtered solvable polynomial algebra with filtration
$FA=\{ F_pA\}_{p\in\mathbb{N}}$, where solvable polynomial algebras
are in the sense of (A.~Kandri-Rody and V.~Weispfenning,
Non-commutative Gr\"obner bases in algebras of solvable type. {\it
J. Symbolic Comput.}, 9(1990), 1--26), and $FA$ is constructed with
respect to a positive-degree function $d(~)$ on $A$. By introducing
minimal F-bases and minimal standard bases respectively for left
$A$-modules and their submodules with respect to good filtrations,
minimal filtered free resolutions for finitely generated $A$-modules
are introduced. It is shown that any two minimal F-bases,
respectively any two minimal standard bases have the same number of
elements and the same number of elements of the same filtered
degree; that minimal filtered free resolutions are unique up to
strict filtered isomorphism of chain complexes in the category of
filtered $A$-modules; and that  minimal finite filtered free
resolutions can be algorithmically computed by employing Gr\"obner
basis theory for modules over $A$ with respect to any graded left
monomial ordering on free left $A$-modules.  }
\end{minipage}\end{center} {\parindent=0pt\vskip 6pt

 \vskip 6pt

}

\def\v5{\vskip .5truecm}\def\QED{\hfill{$\Box$}}\def\hang{\hangindent\parindent}
\def\textindent#1{\indent\llap{#1\enspace}\ignorespaces}
\def\item{\par\hang\textindent}
\def \r{\rightarrow}\def\OV#1{\overline {#1}}
\def\normalbaselines{\baselineskip 24pt\lineskip 4pt\lineskiplimit 4pt}
\def\mapdown#1{\llap{$\vcenter {\hbox {$\scriptstyle #1$}}$}
                                \Bigg\downarrow}
\def\mapdownr#1{\Bigg\downarrow\rlap{$\vcenter{\hbox
                                    {$\scriptstyle #1$}}$}}
\def\mapright#1#2{\smash{\mathop{\longrightarrow}\limits^{#1}_{#2}}}
\def\NZ{\mathbb{N}}

\def\LH{{\bf LH}}\def\LM{{\bf LM}}\def\LT{{\bf
LT}}\def\KX{K\langle X\rangle} \def\KS{K\langle X\rangle}
\def\B{\mathscr{B}} \def\LC{{\bf LC}} \def\G{{\cal G}} \def\FRAC#1#2{\displaystyle{\frac{#1}{#2}}}
\def\SUM^#1_#2{\displaystyle{\sum^{#1}_{#2}}} \def\O{{\cal O}}  \def\J{{\bf J}}\def\BE{\B (e)}
\def\PRCVE{\prec_{\varepsilon\hbox{-}gr}}\def\BV{\B (\varepsilon )}\def\PRCEGR{\prec_{e\hbox{\rm -}gr}}

\def\KS{K\langle X\rangle}
\def\LR{\langle X\rangle}
\def\HL{{\rm LH}}\def\T#1{\widetilde #1}

\vskip 1truecm

\section*{1. Introduction and Preliminaries}
In [Li3] it has been shown that the methods and algorithms,
developed  in ([CDNR], [KR]) for computing minimal homogeneous
generating sets  of graded submodules and graded quotient modules of
free modules  over a commutative polynomial algebra, can be adapted
for computing minimal homogeneous generating sets  of graded
submodules and graded quotient modules of free modules over a
weighted $\mathbb{N}$-graded (noncommutative) solvable polynomial
algebra, and consequently, algorithmic procedures for computing
minimal finite graded free resolutions over weighted
$\mathbb{N}$-graded solvable polynomial algebras are achieved, where
a solvable polynomial algebra $A$ is in the sense of [K-RW] and a
weighted $\NZ$-gradation attached to $A$ is specified to that
constructed with respect to a  positive-degree function $d(~)$ on
the PBW $K$-basis of $A$. Also we learnt from ([KR], Definition
4.2.13, Theorem 4.3.19) that if $L=\oplus_{i=1}^sAe_i$ is a free
module over a commutative polynomial algebra $A=K[x_1,\ldots ,x_n]$,
then a minimal Macaulay basis  for {\it any} finitely generated
submodule $N$ of  $L$  can be obtained via  computing a minimal
homogeneous generating set of the graded submodule generated by all
leading homogeneous elements (degree forms) of $N$ taking with
respect to a fixed gradation for $A$ and a fixed gradation for $L$,
or can be obtained via computing a minimal homogeneous generating
set of the graded submodule generated by all homogenized elements of
$N$ with respect to the variable $x_0$ in $\OV A=K[x_0,x_1,\ldots
,x_n]$ (though Theorem 4.3.19 in [KR] is mentioned only for ideals
of $A$,  this is true for submodules of $L$ as we will see in later
Section 5 where $\OV A$ is replaced by the Rees algebra $\T A$ of
$A$). It is well known that Macaulay bases in more general context
are called {\it standard bases}, especially standard bases for
(two-sided) ideals in general {\it ungraded} commutative and
noncommutative algebras were introduced by E.S. Golod ([Gol], 1986)
in terms of the $\Gamma$-filtered structures attached to arbitrary
associative  algebras, where $\Gamma$ is a well-ordered semigroup
with respect to a well-ordering on $\Gamma$. Motivated by the work
of [Gol],  the results ([KR], Proposition 4.2.15, Proposition
4.3.21, Theorem 4.6.3, Proposition 4.7.24) and the work of [Li3],
the aim  of this paper is first to introduce minimal filtered free
resolutions over a weighted $\NZ$-filtered solvable polynomial
algebra $A$, and then,  to give algorithmic procedures for computing
minimal finite filtered free resolutions in the case where graded
left monomial orderings on free left $A$-modules are used. Our goal
will be reached by employing the filtered-graded transfer techniques
and the Gr\"obner basis theory for solvable polynomial algebras as
well as  their modules. More clearly, the contents of this paper are
arranged as follows:\par

1. Introduction and Preliminaries\par

2. Filtered Free Modules over Weighted $\NZ$-Filtered Solvable
Polynomial Algebras\par

3. Filtered-Graded Transfer of Left Gr\"obner Bases for Modules\par

4. F-Bases and Standard Bases with Respect to  Good Filtrations\par

5. Computation of Minimal F-Bases and Minimal Standard Bases\par

6. The Uniqueness of Minimal Filtered Free Resolutions\par

7. Computation of  Minimal Finite Filtered Free Resolutions\v5

Throughout  this paper, $K$ denotes a field, $K^*=K-\{ 0\}$;
$\mathbb{N}$ denotes the additive monoid of nonnegative integers,
and $\mathbb{Z}$ denotes the additive group of integers; all
algebras are associative $K$-algebras with the multiplicative
identity 1, and modules over an algebra are meant left unitary
modules. \v5

We start by recalling briefly some basics on Gr\"obner basis theory
for solvable polynomial algebras and their modules. The main
references are [AL1], [Gal], [K-RW], [Kr], [LW], [Li1], and [Lev].
Let  $A=K[a_1,\ldots ,a_n]$ be a finitely generated $K$-algebra with
the {\it minimal set of generators} $\{ a_1,\ldots ,a_n\}$. If, for
some permutation $\tau =i_1i_2\cdots i_n$ of $1,2,\ldots ,n$, the
set $\B =\{ a^{\alpha}=a_{i_1}^{\alpha_1}\cdots
a_{i_n}^{\alpha_n}~|~\alpha =(\alpha_1,\ldots ,\alpha_n)\in\NZ^n\}
,$ forms a $K$-basis of $A$, then $\B$ is referred to as a {\it PBW
$K$-basis} of $A$. It is clear that if $A$ has a PBW $K$-basis, then
we can always assume that $i_1=1,\ldots ,i_n=n$. Thus, we make the
following convention once for all.{\parindent=0pt\v5

{\bf Convention} From now on in this paper, if we say that an
algebra $A$ has the PBW $K$-basis $\B$, then it means that
$$\B =\{ a^{\alpha}=a_{1}^{\alpha_1}\cdots
a_{n}^{\alpha_n}~|~\alpha =(\alpha_1,\ldots ,\alpha_n)\in\NZ^n\} .$$
Moreover, adopting the commonly used terminology in computational
algebra, elements of $\B$ are referred to as {\it monomials} of
$A$.} \v5

Suppose that $A$ has the PBW $K$-basis $\B$ as presented above and
that  $\prec$ is a total ordering on $\B$. Then every nonzero
element $f\in A$ has a unique expression
$$f=\lambda_1a^{\alpha (1)}+\lambda_2a^{\alpha (2)}+\cdots +\lambda_ma^{\alpha (m)},~
\lambda_j\in K^*,~a^{\alpha
(j)}=a_1^{\alpha_{1j}}a_2^{\alpha_{2j}}\cdots a_n^{\alpha_{nj}}\in\B
,~1\le j\le m.$$ If  $a^{\alpha (1)}\prec a^{\alpha (2)}\prec\cdots
\prec a^{\alpha (m)}$ in the above representation, then the {\it
leading monomial of $f$} is defined as $\LM (f)=a^{\alpha (m)}$, the
{\it leading coefficient of $F$} is defined as $\LC (f)=\lambda_m$,
and the {\it leading term of $f$} is defined as $\LT
(f)=\lambda_ma^{\alpha (m)}$. {\parindent=0pt\v5

{\bf 1.1. Definition}  Suppose that the $K$-algebra $A=K[a_1,\ldots
,a_n]$ has the PBW $K$-basis $\B$. If $\prec$ is a total ordering on
$\B$ that satisfies the following three
conditions:}{\parindent=1.3truecm\par

\item{(1)} $\prec$ is a well-ordering;\par

\item{(2)} For $a^{\gamma},a^{\alpha},a^{\beta}, a^{\eta}\in\B$, if
$a^{\alpha}\prec a^{\beta}$ and $\LM
(a^{\gamma}a^{\alpha}a^{\eta})$, $\LM
(a^{\gamma}a^{\beta}a^{\eta})\not\in K$, then $\LM
(a^{\gamma}a^{\alpha}a^{\eta})\prec\LM
(a^{\gamma}a^{\beta}a^{\eta})$;\par

\item{(3)} For $a^{\gamma},a^{\alpha},a^{\beta},a^{\eta}\in\B$, if $a^{\beta}\ne
a^{\gamma}$, and $a^{\gamma}=\LM (a^{\alpha}a^{\beta}a^{\eta})$,
then $a^{\beta}\prec a^{\gamma}$ (thereby $1\prec a^{\gamma}$ for
all $a^{\gamma}\ne 1$),\par}{\parindent=0pt

then $\prec$ is called a {\it monomial ordering} on $\B$ (or a
monomial ordering on $A$). }\par

If $\prec$ is a monomial ordering on $\B$, then we call $(\B
,\prec)$ an {\it admissible system} of $A$.\v5

Note that if a $K$-algebra $A=K[a_1,\ldots ,a_n]$ has the PBW
$K$-basis $\B =\{ a^{\alpha}=a_1^{\alpha_1}\cdots
a_n^{\alpha_n}~|~\alpha =(\alpha_1,\ldots ,\alpha_n)\in\NZ^n\}$,
then for any given $n$-tuple $(m_1,\ldots ,m_n)\in\mathbb{N}^n$, a
{\it weighted degree function} $d(~)$ is well defined on nonzero
elements of $A$, namely, for each $a^{\alpha}=a_1^{\alpha_1}\cdots
a_n^{\alpha_n}\in\B$, $d(a^{\alpha})=m_1\alpha_1+\cdots
+m_n\alpha_n$, and for each nonzero
$f=\sum_{i=1}^s\lambda_ia^{\alpha (i)}\in A$ with $\lambda_i\in K^*$
and $a^{\alpha (i)}\in\B$, $d(f)=\max\{ d(a^{\alpha (i)})~|~1\le
i\le s\}.$ If $d(a_i)=m_i>0$ for $1\le i\le n$, then $d(~)$ is
referred to as a {\it positive-degree function} on $A$. \par

Let  $d(~)$ be a  positive-degree function on $A$. If $\prec$ is a
monomial ordering on $\B$ such that for all
$a^{\alpha},a^{\beta}\in\B$,
$$a^{\alpha}\prec a^{\beta}~\hbox{implies}~d(a^{\alpha})\le d(a^{\beta}),\leqno{(*)}$$
then we call $\prec$ a {\it graded monomial ordering} with respect
to $d(~)$, and from now on, unless otherwise stated we always use
$\prec_{gr}$ to denote a  graded monomial ordering.\v5

As one may see from  the literature that in both the commutative and
noncommutative computational algebra, the most popularly used graded
monomial orderings on an algebra $A$ with the PBW $K$-basis $\B$ are
those graded (reverse) lexicographic orderings with respect to the
degree function $d(~)$ such that $d(a_i)=1$, $1\le i\le
n$.{\parindent=0pt\v5

{\bf 1.2. Definition} Suppose that the $K$-algebra $A=K[a_1,\ldots
,a_n]$ has an admissible system $(\B ,\prec )$. If for all
$a^{\alpha}=a_1^{\alpha_1}\cdots a_n^{\alpha_n}$,
$a^{\beta}=a_1^{\beta_1}\cdots a^{\beta_n}_n\in\B$, the following
condition is satisfied:
$$\begin{array}{rcl} a^{\alpha}a^{\beta}&=&\lambda_{\alpha ,\beta}a^{\alpha +\beta}+f_{\alpha ,\beta},\\
&{~}&\hbox{where}~\lambda_{\alpha ,\beta}\in K^*,~a^{\alpha
+\beta}=a_1^{\alpha_1+\beta_1}\cdots a_n^{\alpha_n+\beta_n},~\hbox{and}\\
&{~}&f_{\alpha ,\beta}\in K\hbox{-span}\B~\hbox{with}~\LM (f_{\alpha
,\beta})\prec a^{\alpha +\beta}~\hbox{whenever}~f_{\alpha ,\beta}\ne
0,\end{array}$$ then $A$ is said to be a {\it solvable polynomial
algebra}. \v5

{\bf Remark} Let $A=K[a_1,\ldots ,a_n]$ be a finitely generated
$K$-algebra and $K\langle X\rangle =K\langle X_1,\ldots ,X_n\rangle$
the free $K$-algebra on $\{X_1,\ldots ,X_n\}$. Then it follows from
[Li2] that $A$ is a solvable polynomial algebra if and only
if{\parindent=1.3truecm\par

\item{(1)} $A\cong K\langle X\rangle /\langle G\rangle$ with a finite
set of defining relations $G=\{ g_1,\ldots ,g_m\}$ such that with
respect to some monomial ordering $\prec_{_X}$ on $\KX$, $G$ is a
Gr\"obner basis  of the ideal $\langle G\rangle$, and the set of
normal monomials (mod $G$) gives rise to a PBW $K$-basis $\B$ for
$A$, and\par

\item{(2)} there is a monomial ordering $\prec$ on $\B$ (not necessarily
the restriction of $\prec_{_X}$ on $\B$) such that the condition on
monomials given in  Definition 1.2 is satisfied (see Example (3)
given in the next section for an
illustration).\par}{\parindent=0pt\par

Thus, solvable polynomial algebras are completely determinable and
constructible in a computational way.}}\v5

By Definition 1.2 it is straightforward that if $A$ is a solvable
polynomial algebra and $f,g\in A$ with $\LM (f)=a^{\alpha}$, $\LM
(g)=a^{\beta}$, then
$$\LM (fg)=\LM (\LM (f)\LM (g))=\LM (a^{\alpha}a^{\beta})=a^{\alpha +\beta}.\leqno{(\mathbb{P}1)}$$
We shall freely use this property in the rest of this paper without
additional indication.\par

The results mentioned in the Theorem below are summarized from
([K-RW], Sections 2 -- 5).{\parindent=0pt\v5

{\bf 1.3. Theorem}  Let $A=K[a_1,\ldots ,a_n]$ be a solvable
polynomial algebra with admissible system $(\B ,\prec )$. The
following statements hold.\par

(i) $A$ is a (left and right) Noetherian domain.\par

(ii) With respect to the given $\prec$ on $\B$, every nonzero left
ideal $I$ of $A$ has a finite  left Gr\"obner basis $\G=\{
g_1,\ldots ,g_t\}\subset I$ in the sense
that}{\parindent=.5truecm\par

\item{$\bullet$} if $f\in I$ and $f\ne 0$, then there is a $g_i\in\G$ such that $\LM (g_i)|\LM (f)$,
i.e., there is some $a^{\gamma}\in\B$ such that $\LM (f)=\LM
(a^{\gamma}\LM (g_i))$, or equivalently, with $\gamma
(i_j)=(\gamma_{i_{1j}},\gamma_{i_{2j}},\ldots
,\gamma_{i_{nj}})\in\NZ^n$, $f$ has a left Gr\"obner representation:
$$\begin{array}{rcl} f&=&\sum_{i,j}\lambda_{ij}a^{\gamma (i_j)}g_j,~\hbox{where}~\lambda_{ij}\in K^*,
~a^{\gamma (i_j)}\in\B ,~
g_j\in \G,\\
&{~}&\hbox{satisfying}~\LM (a^{\gamma (i_j)}g_j)\preceq\LM
(f)~\hbox{for all}~(i,j).\end{array}$$}{\parindent=0pt\par

(iii) The Buchberger algorithm, that computes a finite Gr\"obner
basis for a finitely generated commutative polynomial ideal, has a
complete noncommutative version that computes a finite left
Gr\"obner basis for a finitely generated left ideal
$I=\sum_{i=1}^mAf_i$ of $A$ (see {\bf Algorithm 1} given in the end
of this section).\par

(iv) Similar results of (ii) and (iii) hold for right ideals and
two-sided ideals of $A$.\par\QED}\v5

Let $A=K[a_1,\ldots ,a_n]$ be a solvable polynomial algebra with
admissible system $(\B ,\prec )$, and let $L=\oplus_{i=1}^sAe_i$ be
a free (left) $A$-module with the $A$-basis $\{ e_1,\ldots ,e_s\}$.
Then $L$ is a Noetherian module with the $K$-basis $$\BE =\{
a^{\alpha}e_i~|~a^{\alpha}\in\B ,~1\le i\le s\} .$$ We also call
elements of $\BE$ {\it monomials} in $L$. If $\prec_{e}$ is a total 
ordering on $\BE$, and if $\xi =\sum_{j=1}^m\lambda_ja^{\alpha 
(j)}e_{i_j}\in L$, where $\lambda_j\in K^*$ and $\alpha 
(j)=(\alpha_{j_1},\ldots ,\alpha_{j_n})\in\NZ^n$, such that 
$a^{\alpha (1)}e_{i_1}\prec_{e} a^{\alpha 
(2)}e_{i_2}\prec_{e}\cdots\prec_{e} a^{\alpha (m)}e_{i_m}$, then by 
$\LM (\xi )$ we denote the {\it leading monomial} $a^{\alpha 
(m)}e_{i_m}$ of $\xi $, by $\LC (\xi )$ we denote the {\it leading 
coefficient} $\lambda_m$ of $\xi $,  and by $\LT (\xi )$ we denote 
the {\it leading term} $\lambda_ma^{\alpha (m)}e_{i_m}$ of $f$.
\par

With respect to the given monomial ordering $\prec$ on $\B$, a total
ordering $\prec_{e}$ on $\BE$ is called a {\it left monomial 
ordering} if the following two conditions are satisfied:
\par

(1) $a^{\alpha}e_i\prec_{e} a^{\beta}e_j$ implies  $\LM
(a^{\gamma}a^{\alpha}e_i)\prec_{e} \LM (a^{\gamma}a^{\beta}e_j)$ for 
all $a^{\alpha}e_i$, $a^{\beta}e_j\in\BE$, $a^{\gamma}\in\B$;\par

(2) $a^{\beta}\prec a^{\beta}$ implies $a^{\alpha}e_i\prec_{e}  
a^{\beta}e_i$ for all $a^{\alpha},a^{\beta}\in\B$ and $1\le i\le s$. 
{\parindent=0pt\par

From the definition it is straightforward to check that every left
monomial ordering  $\prec_{e}$ on $\BE$ is a well-ordering. 
Moreover, if  $f\in A$ with $\LM (f)=a^{\gamma}$ and $\xi\in L$ with 
$\LM (\xi )=a^{\alpha}e_i$, then by referring to the foregoing 
$(\mathbb{P}1)$ we have
$$\LM (f\xi )=\LM (\LM (f)\LM (\xi ))=\LM (a^{\gamma}a^{\alpha}e_i)=a^{\gamma +\alpha}e_i.\leqno{(\mathbb{P}2)}$$
We shall also freely use this property in the rest of this paper
without additional indication.}\par

Actually as in the commutative case (cf. [AL2], [KR]), any left
monomial ordering $\prec$ on $\B$ may induce  two  left monomial
orderings on $\BE$:
$$\begin{array}{l} (\hbox{{\bf TOP} ordering})\quad a^{\alpha}e_i\prec_{e} a^{\beta}e_j
\Longleftrightarrow a^{\alpha}\prec a^{\beta},~\hbox{or}~
a^{\alpha}=a^{\beta}~
\hbox{and}~i<j;\\
(\hbox{{\bf POT} ordering})\quad a^{\alpha}e_i\prec_{e}
a^{\beta}e_j\Longleftrightarrow i<j,~
\hbox{or}~i=j~\hbox{and}~a^{\alpha}\prec a^{\beta}.\end{array}$$ \v5

Let $\prec_{e}$ be a left monomial ordering  on the $K$-basis $\BE$ 
of $L$, and $a^{\alpha}e_i$, $a^{\beta}e_j\in\BE$, where $\alpha 
=(\alpha_1,\ldots ,\alpha_n)$, $\beta =(\beta_1,\ldots 
,\beta_n)\in\NZ^n$. We say that {\it $a^{\alpha}e_i$ divides 
$a^{\beta}e_j$}, denoted $a^{\alpha}e_i|a^{\beta}e_j$, if $i=j$ and 
$a^{\beta}e_i=\LM (a^{\gamma}a^{\alpha}e_i)$ for some 
$a^{\gamma}\in\B$. It follows from the foregoing property 
$(\mathbb{P}2)$ that
$$a^{\alpha}e_i|a^{\beta}e_j~\hbox{if and only if}~i=j~\hbox{and}~\beta_i\ge
\alpha_i,~1\le i\le n.$$ This division of monomials can be extended
to a division algorithm of dividing an element $\xi$ by a finite
subset of nonzero elements  $U =\{\xi_1,\ldots ,\xi_m\}$ in $L$.
That is, if there is some $\xi_{i_1}\in U$ such that $\LM
(\xi_{i_1})|\LM (\xi )$, i.e., there is a monomial $a^{\alpha
(i_1)}\in\B$ such that $\LM (\xi )=\LM (a^{\alpha (i_1)}\xi_{i_1})$,
then $\xi ':=\xi -\frac{\LC (\xi )}{\LC (a^{\alpha
(i_1)}\xi_{i_1})}a^{\alpha (i_1)}\xi_{i_1}$; otherwise, $\xi ':=\xi
-\LT (\xi )$. Executing this procedure for $\xi '$ and so on, it
follows from the well-ordering property of $\prec_{e}$ that after 
finitely many repetitions  $\xi$ has an expression
$$\begin{array}{rcl} \xi&=&\sum_{i,j}\lambda_{ij}a^{\alpha (i_j)}\xi_j+\eta ,~\hbox{where}~\lambda_{ij}\in K,
~a^{\alpha (i_j)}\in\B ,~\xi_j\in U ,\\
&{~}&\eta =0~\hbox{or}~\eta =\sum_k\lambda_ka^{\gamma
(k)}e_k~\hbox{with}~\lambda_k\in K,~a^{\gamma (k)}e_k\in\B (e)
,\\
&{~}&\hbox{satisfying}\\
&{~}&\LM (a^{\alpha (i_j)}\xi_j)\preceq_{e}\LM (\xi )~\hbox{for
all}~\lambda_{ij}\ne 0,~\hbox{and if}~\eta\ne 0,~\hbox{then}\\
&{~}&a^{\gamma (k)}e_k\preceq_{e}\LM (\xi),~\LM (\xi_i){\not 
|}~a^{\gamma (k)}e_k~\hbox{for all}~\xi_i\in U~\hbox{and 
all}~\lambda_k\ne 0.\end{array}$$ The element $\eta$ appeared in the 
above expression is called a {\it remainder} of $\xi$ on division by 
$ U$, and is usually denoted by $\OV{\xi}^{ U}$, i.e., $\OV{\xi}^{ 
U}=\eta$. If $\OV{\xi}^{ U}=0$, then we say that $\xi$ is {\it 
reduced to zero}  on division by $ U$. A nonzero element $\xi\in L$ 
is said to be {\it normal} (mod $ U$) if $\xi =\OV{\xi}^{ U}$.\par

Based on the division algorithm, the notion of a {\it left Gr\"obner
basis} for a submodule $N$ of the free  module
$L=\oplus_{i=1}^sAe_i$ comes into play. Since $A$ is a Noetherian
domain, it follows that  $L$ is a Noetherian $A$-module and the
following proposition holds. {\parindent=0pt\v5

{\bf 1.4. Theorem} With respect to the given $\prec_{e}$ on $\BE$, 
every nonzero submodule $N$ of $L$ has a finite left Gr\"obner basis 
$\G =\{ g_1,\ldots ,g_m\}\subset N$ in the sense that 
{\parindent=.5truecm\par

\item{$\bullet$} if $\xi\in N$ and $\xi\ne 0$,
then $\LM (g_i)|\LM (\xi )$ for some $g_i\in\G$, i.e., there is a
monomial $a^{\gamma}\in\B$ such that $\LM (\xi )=\LM (a^{\gamma}\LM
(g_i))$, or equivalently, $\xi$ has a {\it left Gr\"obner
representation} $\xi =\sum_{i,j}\lambda_{ij}a^{\alpha (i_j)}g_j$,
where $\lambda_{ij}\in K^*$, $a^{\alpha (i_j)}\in\B$ with $\alpha
(i_j)=(\alpha_{i_{j1}},\ldots ,\alpha_{i_{jn}})\in\NZ^n$,
$g_j\in\G$,  satisfying $\LM (a^{\alpha (i_j)}g_j)\preceq_{e}\LM 
(\xi)$;\par}

moreover, starting with any finite generating set of $N$, such a
left Gr\"obner basis $\G$ can be computed by running a
noncommutative version of the commutative Buchberger algorithm (cf.
[Bu1], [Bu2], [BW]) for modules over solvable polynomial
algebras.\par\QED}\v5

For the the use of later Section 7, we recall the noncommutative
version of the Buchberger algorithm for modules over solvable
polynomial algebras  as follows.\par

Let $N=\sum_{i=1}^mA\xi_i$ with $ U =\{\xi_1,\ldots ,\xi_m\}\subset 
L$.  For  $\xi_i,\xi_j\in U$ with $1\le i<j\le m$, $\LM 
(\xi_i)=a^{\alpha (i)}e_{p}$, $\LM (\xi_j)=a^{\alpha (j)}e_{q}$, 
where $\alpha (i)=(\alpha_{i_1},\ldots,\alpha_{i_n})$, $\alpha 
(j)=(\alpha_{j_1},\ldots ,\alpha_{j_n})$,  put $\gamma 
=(\gamma_1,\ldots ,\gamma_n)$ with $\gamma_k =\max\{ 
\alpha_{i_k},\alpha_{j_k})$. The {\it left S-polynomial} of $\xi_i$ 
and $\xi_j$ is defined as
$$S_{\ell}(\xi_i,\xi_j)=\left\{\begin{array}{ll}
\displaystyle{\frac{1}{\LC (a^{\gamma -\alpha (i)}\xi_i)}}a^{\gamma
-\alpha (i)}\xi_i-\displaystyle{\frac{1}{\LC (a^{\gamma -\alpha
(j)}\xi_j)}}a^{\gamma
-\alpha (j)}\xi_j,&\hbox{if}~ p=q\\
0,&\hbox{if}~p\ne q.\end{array}\right.$$ {\parindent=0pt\vskip
6pt\def\S{{\cal S}}\newpage

\underline{\bf Algorithm
1~~~~~~~~~~~~~~~~~~~~~~~~~~~~~~~~~~~~~~~~~~~~~~~~~~~~~~~~~~~~~~~~~~~~~~~~~~~~~~~~~~~~~~~~~~~~~~~~~~}\par

$\begin{array}{l} \textsc{INPUT}: ~ U =
\{ \xi_1,...,\xi_m\}\\
\textsc{OUTPUT}:~ \G =\{ g_1,...,g_t\},~\hbox{a left Gr\"obner basis
for}~N=\sum_{i=1}^mA\xi_i,~\\
\textsc{INITIALIZATION}:~ m':=m,~\G :=\{g_1=\xi_1,\ldots ,g_{m'}=\xi_m\}  ,\\
\hskip 4.22truecm \S :=\left\{S_{\ell}(g_i,g_j)~\left
|~\begin{array}{l} g_i,g_j\in\G ,~i<j,~\hbox{and for some}~e_t,\\
\LM (g_i)=a^{\alpha}e_t,~\LM
(g_j)=a^{\beta}e_t\end{array}\right.\right\}\end{array}$

$\begin{array}{l}
\textsc{BEGIN}\\
~~~~\textsc{WHILE}~\S\ne \emptyset\\
~~~~~~~~~~\hbox{Choose any}~S_{\ell}(g_i,g_j)\in\S\\
~~~~~~~~~~\S :=\S -\{ S_{\ell}(g_i,g_j)\}\\
~~~~~~~~~~\OV{S_{\ell}(g_i,g_j)}^{\G}=\eta\\
~~~~~~~~~~~~~\textsc{IF}~\eta\ne 0~\hbox{with}~\LM (\eta )=
a^{\rho}e_k~\textsc{THEN}\\
~~~~~~~~~~~~~~~~~m':=m'+1,~g_{m'}:=\eta \\
~~~~~~~~~~~~~~~~~\S :=\S\cup\{ S_{\ell}(g_j,g_{m'})~|~g_j\in\G,~\LM (g_j )=a^{\nu}e_k\}\\
~~~~~~~~~~~~~~~~~\G :=\G\cup\{g_{m'}\},\\
~~~~~~~~~~~~~\textsc{END}\\
~~~~\textsc{END}\\
\textsc{END}\end{array}$\par
\underline{~~~~~~~~~~~~~~~~~~~~~~~~~~~~~~~~~~~~~~~~~~~~~~~~~~~~~~~~~~~~~~~~~~~~~~~~~~~~~~~
~~~~~~~~~~~~~~~~~~~~~~~~~~~~~~~~~~~~~~~~~~~~~~~~~~} \vskip 6pt

One is referred to the up-to-date computer algebra system
\textsc{Singular} [DGPS] for a package implementing {\bf Algorithm
1}. }\v5

\section*{2. Filtered Free Modules over Weighted $\NZ$-Filtered Solvable Polynomial Algebras}
Recall that the $\NZ$-filtered solvable polynomial algebras with the
generators of degree 1 (especially the quadric solvable polynomial
algebras) were studied in ([LW], [Li1]). In this section, by means
of positive-degree functions we introduce more generally the
weighted $\NZ$-filtered solvable polynomial algebras and filtered
free left modules over such algebras. Since the standard bases we
are going to introduce in terms of good filtrations are
generalization of classical Macaulay bases (see a remark given in
Section 4), while a classical Macaulay basis $V$ is  characterized
in  terms of both the leading homogeneous elements (degree forms) of
$V$ and the homogenized elements of $V$ (cf.  [KR], P.38, P.55),
accordingly  both the associated graded algebra (module) and the
Rees algebra (module)  of an $\NZ$-filtered algebra (module) come
into play in our noncommutative filtered context. All notions,
notations and conventions used in Section 1 are maintained.\v5

Let $A=K[a_1,\ldots a_n]$ be a solvable polynomial algebra with
admissible system $(\B ,\prec )$, where $\B =\{
a^{\alpha}=a_1^{\alpha_1}\cdots a_n^{\alpha_n}~|~\alpha
=(\alpha_1,\ldots ,\alpha_n)\in\NZ^n\}$ is the PBW $K$-bsis of $A$
and $\prec$ is a monomial ordering on $\B$, and let $d(~)$ be a
positive-degree function on $A$ such that $d(a_i)=m_i>0$, $1\le i\le
n$ (see Section 2). Put
$$F_pA=K\hbox{-span}\{ a^{\alpha}\in\B~|~d(a^{\alpha})\le p\} ,\quad p\in\NZ,$$
then  it is clear that  $F_pA\subseteq F_{p+1}A$ for all $p\in\NZ$,
$A=\cup_{p\in\NZ}F_pA$, and $1\in F_0A=K$.{\parindent=0pt\v5

{\bf 2.1. Definition} With notation as above, if $F_pAF_qA\subseteq
F_{p+q}A$ holds for all $p,q\in\NZ$, then  we call $A$ a weighted
{\it $\NZ$-filtered solvable polynomial algebra} with respect to the
given positive-degree function $d(~)$ on $A$, and accordingly we
call  $FA=\{ F_pA\}_{p\in\NZ}$ the weighted {\it $\NZ$-filtration}
of $A$. }\v5

Note that a weighted $\NZ$-filtration $FA$ of $A$ is clearly {\it
separated} in the sense that  if $f$ is a nonzero element of $A$,
then $f\in F_pA-F_{p-1}A$ for some $p$. Thus,  if $f\in
F_pA-F_{p-1}A$, then we say that $f$ has {\it filtered degree} $p$
and we use $d_{\rm fil}(f)$ to denote this degree, i.e.,
$$d_{\rm fil}(f)=p\Longleftrightarrow f\in F_pA-F_{p-1}A.\leqno{(\mathbb{P}3)}$$
Since our definition of a weighted $\NZ$-filtration $FA$ for  $A$
depends on a positive-degree function $d(~)$ on $A$, bearing in mind
$(\mathbb{P}3)$ above and the following featured  property of $FA$
will very much help us to deal with the  associated graded
structures of $A$ and filtered $A$-modules.{\parindent=0pt\v5

{\bf 2.2. Lemma}  If $f=\sum_i\lambda_ia^{\alpha (i)}$ with
$\lambda_i\in K^*$ and $a^{\alpha (i)}\in\B$, then $d_{\rm
fil}(f)=p$ if and only if $d(a^{\alpha (i')})=p$ for some $i'$ if
and only if $d(f)=p=d_{\rm fil}(f)$.
\par\QED}\v5

Given a solvable polynomial algebra $A=K[a_1,\ldots ,a_n]$ with
admissible system $(\B ,\prec )$, it follows from Definition 1.2,
Definition 2.1 and Lemma 2.2 that{\parindent=1truecm \vskip6pt

\item{$\bullet$} $A$ is a weighted $\NZ$-filtered solvable polynomial algebra with respect to a given
positive-degree function $d(~)$ if and only if for   $1\le i<j\le
n$, in the relation   $a_ja_i=\lambda_{ji}a_ia_j+f_{ji}$ with
$f_{ji}=\sum\mu_ka^{\alpha (k)}$, $d(a^{\alpha (k)})\le d(a_ia_j)$
holds whenever $\mu_k\ne 0$. \par}{\parindent=0pt\vskip 6pt

This observation helps us to better understand the following
examples.\v5

{\bf Example} (1) If $A=K[a_1,\ldots ,a_n]$ is a weighted
$\NZ$-graded solvable polynomial algebra with respect to a
positive-degre function $d(~)$ on $A$, i.e., $A=\oplus_{p\in\NZ}A_p$
with each $A_p=K$-span$\{a^{\alpha}\in\B~|~d(a^{\alpha})=p\}$ the
degree-$p$ homogeneous part (see [Li3] for more details on weighted
$\NZ$-graded solvable polynomial algebras), then, with respect to
the same positive-degree function $d(~)$ on $A$,  $A$ is turned into
a weighted $\NZ$-filtered solvable polynomial algebra with the
$\NZ$-filtration $FA=\{F_pA\}_{p\in\NZ}$ where each
$F_pA=\oplus_{q\le p}A_q$. \vskip 6pt

{\bf Example} (2) Let $A=K[a_1,\ldots a_n]$ be a solvable polynomial
algebra with the admissible system $(\B ,\prec_{gr})$, where
$\prec_{gr}$ is a graded monomial ordering on $\B$ with respect to a
given positive-degree function $d(~)$ on $A$ (see the definition of
$\prec_{gr}$ given in Section 1). Then by referring to Definition
1.2 and the above observation, one easily sees that $A$ is a
weighted $\NZ$-filtered solvable polynomial algebra with respect to
the same $d(~)$. In the case where $\prec_{gr}$ respects  $d(a_i)=1$
for $1\le i\le n$, Definition 1.2 entails that the generators of $A$
satisfy
$$a_ja_i=\lambda_{ji}a_ia_j+\sum \lambda^{ji}_{k\ell}a_ka_{\ell}+\sum \lambda^{ji}_ta_t+\mu_{ji},~1\le i<j\le n,
~\lambda_{ji}\in K^*,~\lambda^{ji}_{k\ell},\lambda^{ji}_t,
\mu_{ji}\in K.$$ In [Li1] such $\NZ$-filtered solvable polynomial
algebras are referred to as {\it quadric solvable polynomial
algebras} which include numerous significant algebras such as Weyl
algebras and enveloping algebras of Lie algebras. One is referred to
[Li1] for some detailed study of quadric solvable polynomial
algebras  by means of the filtered-graded transfer of Gr\"obner
bases.}\vskip 6pt

The next example provides  weighted $\NZ$-filtered solvable
polynomial algebras in which some generators may have degree $\ge
2$.{\parindent=0pt\v5

{\bf Example} (3)  Considering the $\NZ$-graded structure of the
free $K$-algebra $\KS =K\langle X_1,X_2, X_3\rangle$ by assigning
$X_1$ the degree 2, $X_2$ the degree 1 and $X_3$ the degree 4, let
$I$ be the ideal of $\KS$ generated by the elements
$$\begin{array}{l} g_1=X_1X_2- X_2X_1,\\
g_2=X_3X_1-\lambda X_1X_3-\mu X_3X_2^2-f(X_2),\\
g_3=X_3X_2- X_2X_3,\end{array}$$ where $\lambda\in K^*$, $\mu\in K$,
and  $f(X_2)\in K$-span$\{1,X_2,X_2^2,\ldots ,X_2^6\}$.  If we use
the graded lexicographic ordering
$X_2\prec_{grlex}X_1\prec_{grlex}X_3$ on $\KS$, then it is
straightforward to verify that $\G =\{ g_1,g_2,g_3\}$ forms a
Gr\"obner basis for $I$, and that $\B =\{
a^{\alpha}=a_1^{\alpha_1}a_2^{\alpha_2}a_3^{\alpha_3}~|~\alpha
=(\alpha_1,\alpha_2,\alpha_3)\in\NZ^3\}$ is a PBW basis for the
$K$-algebra $A$, where $A=K[a_1,a_2, a_3]=\KS /I$ with  $a_1$, $a_2$
and $a_3$ denoting the cosets $X_1+I$, $X_2+I$ and $X_3+I$ in $\KS
/I$ respectively.  Since $a_3a_1=\lambda a_1a_3+\mu
a_2^2a_3+f(a_2)$, where $f(a_2)\in K$-span$\{ 1,a_2,a_2^2,\ldots
,a_2^6\}$, we see that $A$ has the monomial ordering $\prec_{lex}$
on $\B$ such that $a_3\prec_{lex}a_2\prec_{lex}a_1$ and $\LM (\mu
a_2^2a_3+f(a_2))\prec_{lex}a_1a_3$, thereby $A$ is turned into a
weighted $\NZ$-filtered solvable polynomial algebra with respect to
$\prec_{lex}$ and the degree function $d(~)$ such that $d(a_1)=2$,
$d(a_2)=1$, and $d(a_3)=4$. Moreover, one may also check that with
respect to the same degree function $d(~)$, the graded lexicographic
ordering $a_3\prec_{grlex}a_2\prec_{grlex}a_1$ is another choice to
make $A$ into a weighted $\NZ$-filtered solvable polynomial
algebra.}\v5

Let $A$ be a weighted  $\NZ$-filtered solvable polynomial algebra
with  admissible system $(\B ,\prec )$, and let
$FA=\{F_pA\}_{p\in\NZ}$ be the $\NZ$-filtration of $A$ constructed
with respect to a given positive-degree function $d(~)$ on $A$. Then
$A$ has the associated $\NZ$-graded $K$-algebra
$G(A)=\oplus_{p\in\NZ}G(A)_p$ with $G(A)_0=F_0A=K$ and
$G(A)_p=F_pA/F_{p-1}A$ for $p\ge 1$, where for $\OV f=f+F_{p-1}A\in
G(A)_p$, $\OV g=g+F_{q-1}A$, the multiplication is given by $\OV
f\OV g=fg+F_{p+q-1}A\in G(A)_{p+q}$. Another $\NZ$-graded
$K$-algebra determined by $FA$ is the Rees algebra $\T A$ of $A$,
which is defined as $\T A=\oplus_{p\in\NZ}\T A_p$ with $\T
A_p=F_pA$, where the multiplication of $\T A$ is induced by
$F_pAF_qA\subseteq F_{p+q}A$, $p, q\in\NZ$. For convenience, we fix
the following notations once for all:{\parindent=.5truecm\par

\item{$\bullet$} If $h\in G(A)_p$ and $h\ne 0$, then
we write $d_{\rm gr}(h)$ for the degree of $h$ as a homogeneous
element of $G(A)$, i.e., $d_{\rm gr}(h)=p$.\par

\item{$\bullet$} If $H\in \T A_p$ and $H\ne 0$, then
we write $d_{\rm gr}(H)$ for the degree of $H$ as a homogeneous
element of $\T A$, i.e., $d_{\rm gr}(H)=p$.\par}

Concerning the $\NZ$-graded structure of $G(A)$, if $f\in A$ with
$d_{\rm fil}(f)=p$, then by Lemma 2.2,  the coset $f+F_{p-1}A$
represented by $f$ in $G(A)_p$ is a nonzero homogeneous element of
degree $p$. If we denote this homogeneous element by $\sigma (f)$
(in the literature it is referred to as the principal symbol of
$f$), then $d_{\rm fil}(f)=p=d_{\rm gr}(\sigma (f))$. However,
considering the Rees algebra $\T A$ of $A$, we note that a nonzero
$f\in F_qA$ represents a homogeneous element of degree $q$ in $\T
A_q$ on one hand, and on the other hand it represents a homogeneous
element of degree $q_1$ in $\T A_{q_1}$, where $q_1=d_{\rm
fil}(f)\le q$. So, for a nonzero $f\in F_pA$, we denote the
corresponding homogeneous element of degree $p$ in $\T A_p$ by
$h_p(f)$, while we use $\T f$ to denote the homogeneous element
represented by $f$ in $\T A_{p_1}$ with $p_1=d_{\rm fil}(f)\le p$.
Thus,  $d_{\rm gr}(\T f)=d_{\rm fil}(f)$, and we see that $h_p(f)=\T
f$ if and only if $d_{\rm fil}(f)=p$.
\par

Furthermore, if we write $Z$ for the homogeneous element of degree 1
in $\T A_1$ represented by the multiplicative identity element 1,
then $Z$ is a central regular element of $\T A$, i.e., $Z$ is not a
divisor of zero and is contained in the center of $\T A$. Bringing
this homogeneous element $Z$ into play,  the homogeneous elements of
$\T A$ are featured as follows:{\parindent=1truecm\vskip6pt

\item{$\bullet$} If $f\in A$ with $d_{\rm fil}(f)=p_1$ then for all $p\ge p_1$,
$h_p(f)=Z^{p-p_1}\T f$. In other words, if $H\in\T A_p$ is any
nonzero homogeneous element of degree $p$, then there is some $f\in
F_pA$ such that $H=Z^{p-d(f)}\T f=\T f+(Z^{p-d(f)}-1)\T
f$.\par}{\parindent=0pt\vskip 6pt

It follows that by sending $H$ to $f+F_{p-1}A$ and sending $H$ to
$f$ respectively, $G(A)\cong\T A/\langle Z\rangle$ as $\NZ$-graded
$K$-algebras and  $A\cong \T A/\langle 1-Z\rangle$ as $K$-algebras
(cf. [AVV], [LVO]).} \v5

Since a solvable polynomial algebra $A$ is necessarily a domain, we
summarize two useful properties concerning the multiplication of
$G(A)$ and $\T A$ respectively  into the following lemma. Notations
are as given before.{\parindent=0pt\v5

{\bf 2.3. Lemma} Let $f,g$ be nonzero elements of $A$ with $d_{\rm
fil}(f)=p_1$, $d_{\rm fil}(g)=p_2$. Then\par

(i)  $\sigma (f)\sigma (g)=\sigma (fg)$;\par

(ii)  $\T f\T g=\widetilde{fg}$. If $p_1+p_2\le p$, then
$h_p(fg)=Z^{p-p_1-p_2}\T f\T g$.\par\QED}\v5

With the preparation made above, the results given in the next
theorem,  which are analogues of those concerning quadric solvable
polynomial algebras in ([LW], Section 3; [Li1], CH.IV), may be
derived in a similar way as in loc. cit., thereby detailed proofs
are omitted. {\parindent=0pt\v5

{\bf 2.4. Theorem}  Let $A=K[a_1,\ldots ,a_n]$ be a solvable
polynomial algebra with the admissible system $(\B ,\prec_{gr})$,
where $\prec_{gr}$ is a graded monomial ordering on $\B$ with
respect to a given positive-degree function $d(~)$ on $A$, thereby
$A$ is a weighted $\NZ$-filtered solvable polynomial algebra with
respect to the same $d(~)$ by the foregoing Example (2), and let
$FA=\{ F_pA\}_{p\in\NZ}$ be the corresponding $\NZ$-filtration of
$A$. Considering the associated graded algebra $G(A)$ as well as the
Rees algebra $\T A$ of $A$, the following statements hold.\par

(i) $G(A)=K[\sigma (a_1),\ldots ,\sigma (a_n)]$,  $G(A)$ has the PBW
$K$-basis
$$\sigma (\B)=\{
\sigma (a)^{\alpha}=\sigma (a_1)^{\alpha_1}\cdots \sigma
(a_n)^{\alpha_n}~|~\alpha =(\alpha_1,\ldots ,\alpha_n)\in\NZ^n\} ,$$
and, by referring to Definition 1.2, for  $\sigma (a)^{\alpha}$,
$\sigma (a)^{\beta}\in\sigma (\B )$ such that
$a^{\alpha}a^{\beta}=\lambda_{\alpha ,\beta}a^{\alpha
+\beta}+f_{\alpha ,\beta}$, where $\lambda_{\alpha ,\beta}\in K^*$,
if $f_{\alpha ,\beta}=0$ then
$$\sigma (a)^{\alpha}\sigma(a)^{\beta}=
\lambda_{\alpha ,\beta}\sigma (a)^{\alpha
+\beta},~\hbox{where}~\sigma (a)^{\alpha +\beta}=\sigma
(a_1)^{\alpha_1+\beta_1}\cdots \sigma (a_n)^{\alpha_n+\beta_n};$$
and in the case where $f_{\alpha ,\beta}=\sum_j\mu^{\alpha
,\beta}_ja^{\alpha (j)}\ne 0$ with $\mu^{\alpha ,\beta}_j\in K$,
$$\sigma (a)^{\alpha}\sigma(a)^{\beta}= \lambda_{\alpha ,\beta}\sigma
(a)^{\alpha +\beta}+\displaystyle{\sum_{d(a^{\alpha
(k)})=d(a^{\alpha +\beta})}}\mu^{\alpha ,\beta}_j\sigma (a)^{\alpha
(k)}.$$ Moreover, the ordering $\prec_{_{G(A)}}$ defined on $\sigma
(\B )$ subject to the rule:
$$\sigma (a)^{\alpha}\prec_{_{G(A)}} \sigma (a)^{\beta}\Longleftrightarrow a^{\alpha}\prec_{gr}a^{\beta},
\quad a^{\alpha},a^{\beta}\in\B,$$  is a graded monomial ordering
with respect to the positive-degree function $d(~)$ on $G(A)$ such
that $d(\sigma (a_i))=d(a_i)$ for $1\le i\le n$, that turns $G(A)$
into a weighted $\NZ$-graded  solvable polynomial algebra.
\par

(ii) $\T A=K[\T a_1,\ldots ,\T a_n, Z]$ where $Z$ is the central
regular element of degree 1 in $\T A_1$ represented by 1,  $\T A$
has the PBW $K$-basis $$\T{\B}=\{\T a^{\alpha}Z^m=\T
a_1^{\alpha_1}\cdots \T a_n^{\alpha_n}Z^m~|~\alpha =(\alpha_1,\ldots
,\alpha_n )\in\NZ^n,m\in\NZ\} ,$$ and, by referring to Definition
1.2, for  $\sigma (a)^{\alpha}$, $\sigma (a)^{\beta}\in\sigma (\B )$
such that $a^{\alpha}a^{\beta}=\lambda_{\alpha ,\beta}a^{\alpha
+\beta}+f_{\alpha ,\beta}$, where $\lambda_{\alpha ,\beta}\in K^*$,
if $f_{\alpha ,\beta}=0$ then
$$\T a^{\alpha}Z^s\cdot\T a^{\beta}Z^t=
\lambda_{\alpha ,\beta}\T a^{\alpha +\beta}Z^{s+t},~\hbox{where}~\T
a^{\alpha +\beta}=\T a_1^{\alpha_1+\beta_1}\cdots \T
a_n^{\alpha_n+\beta_n};$$ and in the case where $f_{\alpha
,\beta}=\sum_j\mu^{\alpha ,\beta}_ja^{\alpha (j)}\ne 0$ with
$\mu^{\alpha ,\beta}_j\in K$,
$$\T a^{\alpha}Z^s\cdot\T a^{\beta}Z^t=
\lambda_{\alpha ,\beta}\T a^{\alpha +\beta}Z^{s+t}+\sum_j\mu^{\alpha
,\beta}_j\T a^{\alpha (j)}Z^{q-m_j},~\hbox{where}~q=d(a^{\alpha
+\beta})+s+t,~m_j=d(a^{\alpha (j)}).$$ Moreover, the ordering
$\prec_{_{\T A}}$ defined on $\T{\B}$ subject to the rule:
$$\T a^{\alpha}Z^s\prec_{_{\T A}}\T a^{\beta}Z^t\Longleftrightarrow a^{\alpha}\prec_{gr}a^{\beta},
\hbox{or}~a^{\alpha}=a^{\beta}~\hbox{and}~s<t,\quad
a^{\alpha},a^{\beta}\in\B,$$ is a monomial ordering on $\T{\B}$
(which is not necessarily a graded monomial ordering), that turns
$\T A$ into a weighted $\NZ$-graded solvable polynomial algebra with
respect to the positive-degree function $d(~)$ on $\T A$ such that
$d(Z)=1$ and $d(\T{a_i})=d (a_i)$ for $1\le i\le n$.\par\QED }\v5

By referring to Lemma 2.2 and  Lemma 2.3, the corollary presented
below is straightforward and will be very often used in discussing
left Gr\"obner bases and standard bases in free left $A$-modules and
their associated graded free $G(A)$-modules as well the graded free
$\T A$-modules (Section 3, Section 4). {\parindent=0pt \v5

{\bf 2.5. Corollary} With the assumption and notations as in Theorem
2.4, if $f=\lambda a^{\alpha}+\sum_j\mu_ja^{\alpha (j)}$ with
$d(f)=p$ and $\LM (f)=a^{\alpha}$, then $p=d_{\rm fil}(f)=d_{\rm
gr}(\sigma (f))=d_{\rm gr}(\T f)$, and
$$\begin{array}{l} \sigma (f)=\lambda\sigma (a)^{\alpha}+\sum_{d(a^{\alpha (j_k)})=p}\mu_{j_k}\sigma (a)^{\alpha (j_k)};\\
\LM(\sigma (f))=\sigma (a)^{\alpha}=\sigma (\LM (f));\\
\T f=\lambda\T a^{\alpha}+\sum_j\mu_j\T a^{\alpha
(j)}Z^{p-d(a^{\alpha (j)})};\\
\LM (\T f)=\T a^{\alpha}=\widetilde{\LM (f)},\end{array}$$ where
$\LM (f)$, $\LM (\sigma (f))$ and $\LM (\T f)$ are taken with
respect to $\prec_{gr}$, $\prec_{_{G(A)}}$ and $\prec_{_{\T A}}$
respectively.\par\QED}\v5

Let $A$ be a weighted $\NZ$-filtered solvable polynomial algebra
with admissible system $(\B ,\prec )$, and let
$FA=\{F_pA\}_{p\in\NZ}$ be the $\NZ$-filtration of $A$ constructed
with respect to a given positive-degree function $d(~)$ on $A$ (see
Section 2). Consider a free $A$-module $L=\oplus_{i=1}^sAe_i$ with
the $A$-basis $\{ e_1,\ldots ,e_s\}$. Then $L$ has the $K$-basis
$\BE =\{ a^{\alpha}e_i~|~a^{\alpha}\in\B ,~1\le i\le s\}$. If  $\{
b_1,\ldots ,b_s\}$ is an {\it arbitrarily} fixed subset of $\NZ$,
then, with $FL=\{ F_qL\}_{q\in\NZ}$ defined by putting
$$F_qL=\{ 0\}~\hbox{if}~q<\min\{ b_1,\ldots ,b_s\};~\hbox{otherwise}~F_qL=
\sum^s_{i=1}\left (\sum_{p_i+b_i\le q}F_{p_i}A\right )e_i,$$ or
alternatively, for $q\ge\min\{ b_1,\ldots, b_s\}$,
$$F_qL=K\hbox{-span}\{
a^{\alpha}e_i\in\BE~|~d(a^{\alpha})+b_i\le q\},$$ $L$ forms an
$\NZ$-{\it filtered free $A$-module} with respect to the
$\NZ$-filtered structure of $A$, that is, every $F_qL$ is a
$K$-subspace of $L$, $F_qL\subseteq F_{q+1}L$ for all $q\in\NZ$,
$L=\cup_{q\in\NZ}F_qL$, $F_pAF_{q}L\subseteq F_{p+q}L$ for all
$p,q\in\NZ$, and for each $i=1,\ldots ,s$,
$$e_i\in F_0L~\hbox{if}~b_i=0;~\hbox{otherwise}~
e_i\in F_{b_i}L-F_{b_i-1}L.$$ {\parindent=0pt\par

{\bf Convention} Let $A$ be a weighted $\NZ$-filtered solvable
polynomial algebra. Unless otherwise stated, from now on in the
subsequent sections  if we say that $L=\oplus_{i=1}^sAe_i$ is a
filtered free $A$-module with the filtration $FL=\{
F_qL\}_{q\in\NZ}$, then $FL$ is always meant the type as constructed
above.}\v5

Let $L=\oplus_{i=1}^sAe_i$ be a filtered free $A$-module with the
filtration $FL=\{ F_qL\}_{q\in\NZ}$, which is constructed with
respect to a given subset $\{ b_1,\ldots ,b_s\}\subset\NZ$. Then
$FL$ is {\it separated} in the sense that if $\xi$ is a nonzero
element of $L$, then $\xi\in F_qL-F_{q-1}L$ for some $q$. Thus, to
make the discussion on $FL$ compatible  with $FA$, if $\xi\in
F_qL-F_{q-1}L$, then we say that $\xi$ has {\it filtered degree} $q$
and we use $d_{\rm fil}(\xi )$ to denote this degree, i.e.,
$$d_{\rm fil}(\xi)=q\Longleftrightarrow \xi\in F_qL-F_{q-1}L.\leqno{(\mathbb{P}4)}$$
For instance, we have $d_{\rm fil}(e_i)=b_i$, $1\le i\le s$.
Comparing with Lemma 2.2 we first note the
following{\parindent=0pt\v5

{\bf 2.6. Lemma} Let $\xi\in L$. Then $d_{\rm fil}(\xi )=q$ if and
only if $\xi =\sum_{i,j}\lambda_{ij}a^{\alpha (i_j)}e_j$, where
$\lambda_{ij}\in K^*$ and $a^{\alpha (i_j)}\in\B$ with $\alpha
(i_j)=(\alpha_{i_{j1}},\ldots ,\alpha_{i_{jn}})\in\NZ^n$, in which
some monomial $a^{\alpha (i_j)}e_j$ satisfy $d(a^{\alpha
(i_j)})+b_j=q$.\par\QED}\v5

Let $L=\oplus_{i=1}^sAe_i$ be a filtered free $A$-module with the
filtration $FL=\{ F_qL\}_{q\in\NZ}$ such that $d_{\rm
fil}(e_i)=b_i$, $1\le i\le s$. Considering the the associated
$\NZ$-graded algebra $G(A)$ of $A$, the filtered free $A$ module $L$
has  the {\it associated $\NZ$-graded $G(A)$-module}
$G(L)=\oplus_{q\in\NZ}G(L)_q$ with $G(L)_q=F_qL/F_{q-1}L$, where for
$\OV f=f+F_{p-1}A\in G(A)_p$, $\OV \xi=\xi+F_{q-1}L\in G(L)_q$, the
module action is given by $\OV f\cdot \OV \xi =f\xi+F_{p+q-1}L\in
G(L)_{p+q}$. As with homogeneous elements in $G(A)$, if $h\in
G(L)_q$ and $h\ne 0$, then  we write $d_{\rm gr}(h)$ for the degree
of $h$ as a homogeneous element of $G(L)$, i.e., $d_{\rm gr}(h)=q$.
If $\xi\in L$ with $d_{\rm fil}(\xi )=q$, then the coset
$\xi+F_{q-1}L$  represented by $\xi$ in $G(L)_q$ is a nonzero
homogeneous element of degree $q$, and if we denote this homogeneous
element by $\sigma (\xi )$ (in the literature it is referred to as
the principal symbol of $\xi$) then $d_{\rm gr}(\sigma (\xi
))=q=d_{\rm fil}(\xi )$.\par

Furthermore, considering the Rees algebra $\T A$ of $A$, the
filtration $FL=\{ F_qL\}_{q\in\NZ}$ of $L$ also defines the {\it
Rees module} $\T L$ of $L$, which is the $\NZ$-graded $\T A$-module
$\T L=\oplus_{q\in\NZ}\T L_q$, where $\T L_q=F_qL$ and the module
action is induced by $F_pAF_qL\subseteq F_{p+q}L$. As with
homogeneous elements in $\T A$, if $H\in \T L_q$ and $H\ne 0$, then
we write $d_{\rm gr}(H)$ for the degree of $H$ as a homogeneous
element of $\T L$, i.e., $d_{\rm gr}(H)=q$. Note that any nonzero
$\xi\in F_qL$ represents a homogeneous element of degree $q$ in $\T
L_q$ on one hand, and on the other hand it represents a homogeneous
element of degree $q_1$ in $\T L_{q_1}$, where $q_1=d_{\rm fil}(\xi
)\le q$. So, for  a nonzero $\xi\in F_qL$ we  denote the
corresponding homogeneous element of degree $q$ in $\T L_q$ by
$h_q(\xi )$,  while we use $\T{\xi}$ to denote the homogeneous
element  represented by $\xi$ in $\T L_{q_1}$ with $q_1=d_{\rm
fil}(\xi )\le q$. Thus,  $d_{\rm gr}(\T{\xi})=d_{\rm fil}(\xi)$, and
we see that $h_q(\xi )=\T{\xi}$ if and only if $d_{\rm fil}(\xi)=q$.
\par

We also note that if  $Z$ denotes the homogeneous element of degree
1 in $\T A_1$ represented by the multiplicative identity element 1,
then, similar to the discussion given before Theorem 2.4, there are
$A$-module isomorphism  $L\cong \T L/(1-Z)\T L$ and graded
$G(A)$-module isomorphism $G(L)\cong \T L/Z\T L$ (cf. [LVVO],
[LVO]). {\parindent=0pt\v5

{\bf 2.7. Lemma} With notation as above, the following statements
hold.\par

(i)  $d_{\rm fil}(f\xi )=d(f)+d_{\rm fil}(\xi )$ holds for all
nonzero $f\in A$ and nonzero $\xi\in L$.\par

(ii) $\sigma (f)\sigma (\xi )=\sigma (f\xi )$ holds for all nonzero
$f\in A$ and nonzero $\xi\in L$.\par

(iii) If $\xi\in L$ with $d_{\rm fil}(\xi )=q\le \ell$, then
$h_{\ell}(\xi )=Z^{\ell -q}\T \xi$. Furthermore, let $f\in A$ with
$d_{\rm fil}(f)=p$, $\xi\in L$ with $d_{\rm fil}(\xi )=q$. Then $\T
f\T \xi =\widetilde{f\xi}$; if $p+q\le \ell$, then $h_{\ell}(f\xi
)=Z^{\ell -p-q}\T f\T \xi$. \vskip 6pt

{\bf Proof} Since $A$ is a solvable polynomial algebra,  $G(A)$ and
$\T A$ are $\NZ$-graded solvable polynomial algebras by Theorem 2.4,
thereby they are necessarily domains.  By the foregoing
$(\mathbb{P}3)$, $(\mathbb{P}4)$ and Lemma 2.6, the verification of
(i) -- (iii) are then straightforward.\QED\v5

{\bf 2.8. Proposition} With notation as fixed before, let
$L=\oplus_{i=1}^sAe_i$ be a filtered free $A$-module with the
filtration $FL=\{ F_qL\}_{q\in\NZ}$ such that $d_{\rm
fil}(e_i)=b_i$, $1\le i\le s$. The following two statements
hold.\par

(i) $G(L)$ is an $\NZ$-graded free $G(A)$-module with the
homogeneous $G(A)$-basis $\{ \sigma (e_1),\ldots ,\sigma (e_s)\}$,
that is,  $G(L)=\oplus_{i=1}^sG(A)\sigma
(e_i)=\oplus_{q\in\NZ}G(L)_q$ with
$$G(L)_q=\sum_{p_i+b_i=q}G(A)_{p_i}\sigma (e_i)\quad q\in\NZ.
$$
Moreover, $\sigma (\B (e))=\{\sigma (a^{\alpha}e_i)=\sigma
(a)^{\alpha}\sigma (e_i)~|~a^{\alpha}e_i\in\B (e)\}$ forms a 
$K$-basis for $G(L)$.\par

(ii) $\T L$ is an $\NZ$-graded free $\T A$-module with the
homogeneous  $\T A$-basis $\{ \T e_1,\ldots ,\T e_s\}$, that is,
$\T L=\oplus_{i=1}^s\T A\T e_i=\oplus_{q\in\NZ}\T L_q$ with
$$\T L_q=\sum_{p_i+b_i=q}\T A_{p_i}\T e_i, \quad q\in\NZ .$$
Moreover, $\widetilde{\B (e)}=\{\T{a}^{\alpha}Z^m\T e_i~|~\T 
a^{\alpha}Z^m\in\T{\B},~1\le i\le s\}$ forms a $K$-basis for $\T L$, 
where $\T{\B}$ is the PBW $K$-basis of $\T A$ determined in Theorem 
4.4(ii).\vskip 6pt

{\bf Proof} Since $d_{\rm fil}(e_i)=b_i$, $1\le i\le s$, if  $\xi
=\sum^s_{i=1}f_ie_i\in F_qL=\sum^s_{i=1}\left (\sum_{p_i+b_i\le
q}F_{p_i}A\right )e_i$, then $d_{\rm fil}(\xi )\le q$. By Lemma 2.7,
$$\begin{array}{l} \sigma (\xi )=\sum_{d(f_i)+b_i=q}\sigma (f_i)\sigma (e_i)
\in\sum_{i=1}^sG(A)_{q-b_i}\sigma (e_i)\\
\\
h_q(\xi )=\sum^s_{i=1}Z^{q-d(f_i)-b_i}\T f_i\T e_i\in\sum^s_{i=1}\T
A_{q-b_i}\T e_i.\end{array}$$ This shows that $\{ \sigma
(e_1),\ldots ,\sigma (e_s)\}$ and $\{ \T e_1,\ldots ,\T e_s\}$
generate the $G(A)$-module $G(L)$ and the $\T A$-module $\T L$,
respectively. Next, since each $\sigma (e_i)$ is a homogeneous
element of degree $b_i$, if a degree-$q$ homogeneous element
$\sum_{i=1}^s\sigma (f_i)\sigma (e_i)=0$, where $f_i\in A$, $d_{\rm
fil}(f_i)+b_i=q$, $1\le i\le s$, then $\sum_{i=1}^sf_ie_i\in
F_{q-1}L$ and hence each $f_i\in F_{q-1-b_i}A$ by Lemma 2.6, a
contradiction. It follows that $\{ \sigma (e_1),\ldots ,\sigma
(e_s)\}$ is linearly independent over $G(A)$. Concerning the linear
independence of $\{ \T e_1,\ldots ,\T e_s\}$ over $\T A$, since each
$\T e_i$ is a homogeneous element of degree $b_i$, if a degree-$q$
homogeneous element $\sum_{i=1}^sh_{p_i}(f_i)\T e_i=0$, where
$f_i\in F_{p_i}A$ and $p_i+b_i=q$, $1\le i\le s$, then
$\sum_{i=1}^sf_ie_i=0$ in $F_qL$ and consequently all $f_i=0$,
thereby $h_{p_i}(f_i)=0$ as desired. Finally, if $\xi\in F_qL$ with
$d_{\rm fil}(\xi )=q$, then by Lemma 2.6, $\xi
=\sum_{i,j}\lambda_{ij}a^{\alpha (i_j)}e_j$ with $\lambda_{ij}\in
K^*$ and $d(a^{\alpha (i_j)})+b_j=\ell_{ij}\le q$. It follows from
Lemma 2.7 that
$$\begin{array}{l} \sigma (\xi )=\sum_{\ell_{ik}=q}
\lambda_{ik}\sigma (a)^{\alpha (i_k)}\sigma (e_k),\\
\T{\xi}=\sum_{i,j}\lambda_{ij}Z^{q-\ell_{ij}}\T a^{\alpha (i_j)}\T
e_j.\end{array}$$ Therefore, a further application of Lemma 2.6 and
Lemma 2.7 shows that $\sigma (\B (e))$ and $\T{\BE}$ are $K$-bases 
for $G(L)$ and $\T L$ respectively.} \v5

\section*{3. Filtered-Graded Transfer of Left Gr\"obner Bases for Modules}
\def\PRCEGR{\prec_{e\hbox{\rm -}gr}}
Throughout this section, we let $A=K[a_1,\ldots ,a_n]$ be a solvable
polynomial algebra with the admissible system $(\B ,\prec_{gr})$,
where $\B=\{ a^{\alpha}=a_1^{\alpha_1}\cdots a_n^{\alpha_n}~|~\alpha
=(\alpha_1,\ldots ,\alpha_n)\in\NZ^n\}$ is  the PBW $K$-basis of $A$
and $\prec_{gr}$ is a graded monomial ordering  with respect to some
given positive-degree function $d(~)$ on $A$ (see Section  2).
Thereby $A$ is turned into a weighted $\NZ$-filtered solvable
polynomial algebra with  the filtration $FA=\{ F_pA\}_{p\in\NZ}$
constructed with respect to the same $d(~)$ (see Example (2) of
Section 2). In order to compute minimal standard bases by employing
both inhomogeneous and homogenous left Gr\"obner bases in later
Section 5, our aim of the current section is to show the relations
between left Gr\"obner bases in a filtered free (left)  $A$-module
$L$ and homogeneous left Gr\"obner bases in $G(L)$ as well as
homogeneous left Gr\"obner bases in $\T L$, which are just module
theory analogues of the results on filtered-graded transfer of
Gr\"obner bases given in ([LW], [Li1]).  All notions, notations and
conventions introduced in previous sections  are maintained.\v5

Let $L=\oplus_{i=1}^sAe_i$ be a filtered free  $A$-module with the
filtration $FL=\{ F_qL\}_{q\in\NZ}$ such that $d_{\rm
fil}(e_i)=b_i$, $1\le i\le s$. Bearing in mind Lemma 2.6, we say
that a left monomial ordering on $\BE$ is a {\it graded left
monomial ordering}, denoted by $\PRCEGR$, if for
$a^{\alpha}e_i,a^{\beta}e_j\in\BE$,
$$a^{\alpha}e_i\PRCEGR a^{\beta}e_j~\hbox{implies}~d_{\rm fil}(a^{\alpha}e_i)=
d(a^{\alpha})+b_i\le d(a^{\beta}) +b_j=d_{\rm fil}(a^{\beta}e_j).$$ 
For instance,  with respect to the given graded monomial ordering 
$\prec_{gr}$ on $\B$ and  the $\NZ$-filtration $FA$ of $A$,   if $\{ 
f_1,\ldots ,f_s\}\subset A$ is a finite subset such that 
$d(f_i)=b_i=d_{\rm fil}(e_i)$, $1\le i\le s$, then, by mimicking the 
Schreyer ordering in the commutative case (see [Sch], or [AL2], 
P.166), one may directly check that the ordering  
$\prec_{s\hbox{-}gr}$ on $\BE$ induced by $\{ f_1,\ldots ,f_s\}$ 
defined subject to the rule: for $a^{\alpha}e_i,a^{\beta}e_j\in\BE$, 
$$a^{\alpha}e_i\prec_{s\hbox{-}gr} a^{\beta}e_j\Longleftrightarrow\left\{\begin{array}{l}
\LM (a^{\alpha}f_i)\prec_{gr}\LM (a^{\beta}f_j),\\
\hbox{or}\\
\LM (a^{\alpha}f_i)=\LM (a^{\beta}f_j)~\hbox{and}~i<j.\\
\end{array}\right.$$
is a graded left monomial ordering on $\BE$. \par

More generally, let $\{ \xi_1,\ldots ,\xi_m\}\subset L$ be a finite
subset, where $d_{\rm fil}(\xi_i)=q_i$, $1\le i\le m$, and let
$L_1=\oplus_{i=1}^mA\varepsilon_i$ be the filtered free $A$-module
with the filtration $FL_1=\{ F_qL_1\}_{q\in \NZ}$ such that  $d_{\rm
fil}(\varepsilon_i)=q_i$, $1\le i\le m$.  Then, given {\it any}
graded left monomial ordering $\PRCEGR$ on $\BE$, the Schreyer
ordering $\prec_{s\hbox{-}gr}$ defined on the $K$-basis $\B
(\varepsilon )=\{ a^{\alpha}\varepsilon_i~|~a^{\alpha}\in\B ,~1\le
i\le m\}$ of $L_1$ subject to the rule: for
$a^{\alpha}\varepsilon_i,a^{\beta}\varepsilon_j\in\B (\varepsilon 
)$,
$$a^{\alpha}\varepsilon_i\prec_{s\hbox{-}gr} a^{\beta}\varepsilon_j
\Longleftrightarrow\left\{\begin{array}{l}
\LM (a^{\alpha}\xi_i)\PRCEGR\LM (a^{\beta}\xi_j),\\
\hbox{or}\\
\LM (a^{\alpha}\xi_i)=\LM (a^{\beta}\xi_j)~\hbox{and}~i<j,\\
\end{array}\right.$$
is a graded left monomial ordering on $\B (\varepsilon )$.\par

Comparing with Lemma 2.2 and Lemma 2.6, the lemma given below
reveals the intrinsic property of a graded left monomial ordering
employed by a filtered free $A$-module.{\parindent=0pt\v5

\def\PRCVE{\prec_{\varepsilon\hbox{-}gr}}\def\BV{\B (\varepsilon )}

{\bf 3.1. Lemma} Let $L=\oplus_{i=1}^sAe_i$ be a filtered free
$A$-module with the filtration $FL=\{ F_qL\}_{q\in\NZ}$ such that
$d_{\rm fil}(e_i)=b_i$, $1\le i\le s$, and let $\prec_{e\hbox{-}gr}$ 
be a graded left monomial ordering on $\BE$. Then $\PRCEGR$ is 
compatible with the filtration $FL$ of $L$ in the sense that $\xi\in 
F_qL-F_{q-1}L$, i.e. $d_{\rm fil}(\xi )=q$, if and only if $\LM (\xi 
)=a^{\alpha}e_i$ with $d_{\rm 
fil}(a^{\alpha}e_i)=d(a^{\alpha})+b_i=q$. \vskip 6pt

{\bf Proof} Let $\xi =\sum_{i,j}\lambda_{ij}a^{\alpha (i_j)}e_j\in
F_qL-F_{q-1}L$. Then  by Lemma 2.6, there is some $a^{\alpha
(i_{\ell})}e_{\ell}$ such that $d(a^{\alpha
(i_{\ell})})+b_{\ell}=q$. If $\LM (\xi )=a^{\alpha (i_t)}e_t$ with
respect to $\PRCEGR$, then $a^{\alpha (i_k)}e_k\PRCEGR a^{\alpha
(i_t)}e_t$ for all $a^{\alpha (i_k)}e_k$ with $k\ne t$. If $\ell
=t$, then $d(a^{\alpha (i_t)})+b_t=q$; otherwise, since $\PRCEGR$ is
a graded left monomial ordering, we have $d(a^{\alpha (i_k)})+b_k\le
d(a^{\alpha (i_t)})+b_t$, in particular, $q=d(a^{\alpha
(i_{\ell})})+b_{\ell}\le d(a^{\alpha (i_t)})+b_t\le q$. Hence
$d_{\rm fil}(a^{\alpha (i_t)}e_t)=d(a^{\alpha (i_t)})+b_t=q$, as
desired.}\par

Conversely, for $\xi =\sum_{i,j}\lambda_{ij}a^{\alpha (i_j)}e_j\in
L$, if,  with respect to $\PRCEGR$, $\LM (\xi )=a^{\alpha (i_t)}e_t$
with $d_{\rm fil}(a^{\alpha (i_t)}e_t)=d(a^{\alpha (i_t) })+b_t=q$,
then $a^{\alpha (i_k)}e_k\PRCEGR a^{\alpha (i_t)}e_t$ for all $k\ne
t$. Since $\PRCEGR$ is a graded left monomial ordering, we have
$d(a^{\alpha (i_k)})+b_k\le d(a^{\alpha (i_t)})+b_t=q$. It follows
from Lemma 2.6 that $d_{\rm fil}(\xi )=q$, i.e.,  $\xi\in
F_qL-F_{q-1}L$. \QED \v5

Let  $L=\oplus_{i=1}^sAe_i$ be a filtered free $A$-module with the
filtration $FL=\{ F_qL\}_{q\in\NZ}$ such that $d_{\rm
fil}(e_i)=b_i$, $1\le i\le s$.  Then, by Proposition 2.8 we know
that the associated graded $G(A)$-module $G(L)$ of $L$ is an
$\NZ$-graded free module, i.e., $G(L)=\oplus_{i=1}^sG(A)\sigma
(e_i)$ with the homogeneous $G(A)$-basis $\{ \sigma (e_1),\ldots
,\sigma (e_s)\}$, and that $G(L)$ has the  $K$-basis $\sigma(\B
(e))=\{ \sigma (a^{\alpha}e_i)=\sigma (a)^{\alpha}\sigma
(e_i)~|~a^{\alpha}e_i\in\B (e)\}$. Furthermore, let
$\prec_{e\hbox{-}gr}$ be a graded left monomial ordering on $\BE$ as 
defined in the beginning of this section. Then we may define an 
ordering $\prec_{\mathbbm{\sigma (e)}\hbox{-}gr}$ on $\sigma(\B 
(e))$ subject to the rule:
$$\sigma (a)^{\alpha}\sigma (e_i)\prec_{\mathbbm{\sigma (e)}\hbox{-}gr}\sigma (a)^{\beta}\sigma (e_j)
\Longleftrightarrow
a^{\alpha}e_i\prec_{e\hbox{-}gr}a^{\beta}e_j,\quad
a^{\alpha}e_i,a^{\beta}e_j\in\B (e).$$ {\parindent=0pt\par
\def\PRCEGR{\prec_{e\hbox{-}gr}}

{\bf 3.2. Lemma} With the ordering $\prec_{\mathbbm{\sigma
(e)}\hbox{-}gr}$ defined above, the following statements hold.\par

(i) $\prec_{\mathbbm{\sigma (e)}\hbox{-}gr}$ is a graded left
monomial ordering on $\sigma(\B (e))$.\par

(ii) (Compare with Corollary 2.5.) $\LM (\sigma (\xi ))=\sigma (\LM
(\xi ))$ holds for all nonzero $\xi\in L$, , where the monomial
orderings used for $\LM (\sigma (\xi ))$ and $\LM (\xi )$ are
$\prec_{\mathbbm{\sigma (e)}\hbox{-}gr}$ and $\PRCEGR$ respectively.
\vskip 6pt

{\bf Proof} (i) Noticing that the given monomial ordering
$\prec_{gr}$ on $A$ is a graded monomial ordering with respect to a
positive-degree function $d(~)$ on $A$, it follows from  Theorem
2.4(i) that  $G(A)$ is turned into an $\NZ$-graded solvable
polynomial algebra by using the graded monomial ordering
$\prec_{_{G(A)}}$ defined on $\sigma (\B )$ subject to the rule:
$\sigma (a)^{\alpha}\prec_{_{G(A)}}\sigma (a)^{\beta}$
$\Longleftrightarrow$ $a^{\alpha}\prec_{gr}a^{\beta}$, where the
positive-degree function on $G(A)$ is given by $d(\sigma
(a_i))=d(a_i)$, $1\le i\le n$. Moreover, since $\sigma (e_i)$ is a
homogeneous element of degree $b_i$ in $G(L)$, $1\le i\le s$, by
Lemma 2.7, it is then straightforward to verify that
$\prec_{\mathbbm{\sigma (e)}\hbox{-}gr}$ is a graded left monomial
ordering on $\sigma(\B (e))$.}\par

(ii) Let $\xi =\sum_{i,j}\lambda_{ij}a^{\alpha (i_j)}e_j$, where
$\lambda_{ij}\in K^*$ and $a^{\alpha (i_j)}\in\B$ with $\alpha
(i_j)=(\alpha_{i_{j1}},\ldots ,\alpha_{i_{jn}})\in\NZ^n$. If $d_{\rm
fil}(\xi )=q$, i.e., $\xi\in F_qL-F_{q-1}L$, then by Lemma 3.1, $\LM
(\xi )=a^{\alpha (i_t)}e_t$ for some $t$ such that $d_{\rm
fil}(a^{\alpha (i_t)}e_t)=d(a^{\alpha (i_t)})+b_t=q$. Since
$\prec_{e\hbox{-}gr}$ is a left graded monomial ordering on $\BE$, 
by Lemma 2.7 we have  $\sigma (\xi )=\lambda_{it}\sigma (a)^{\alpha 
(i_t)}\sigma (e_{t})+\sum_{d(a^{\alpha 
(i_k)})+b_k=q}\lambda_{ik}\sigma (a)^{\alpha (i_k)}\sigma (e_k)$. It 
follows from the definition of $\prec_{\mathbbm{\sigma 
(e)}\hbox{-}gr}$ that $\LM (\sigma (\xi ))=\sigma (a)^{\alpha 
(i_t)}\sigma (e_t)=\sigma (\LM (\xi ))$, as 
desired.\QED{\parindent=0pt\v5

{\bf 3.3. Theorem} Let $N$ be a submodule of the filtered free
$A$-module $L=\oplus_{i=1}^sAe_i$, where $L$ is equipped with the
filtration $FL=\{ F_qL\}_{q\in\NZ}$ such that $d_{\rm
fil}(e_i)=b_i$, $1\le i\le s$, and let $\prec_{e\hbox{-}gr}$ be a 
graded left monomial ordering on $\BE$. For a subset $\G =\{ 
g_1,\ldots ,g_m\}$ of $N$, the following two statements are 
equivalent.\par

(i) $\G$ is a left Gr\"obner basis of $N$ with respect to
$\PRCEGR$.\par

(ii) Putting $\sigma (\G )=\{ \sigma (g_1),\ldots ,\sigma (g_m)\}$
and considering the filtration $FN=\{ F_qN=F_qL\cap N\}_{q\in\NZ}$
of $N$ induced by $FL$ (see Section 2), $\sigma (\G )$ is a left
Gr\"obner basis for the associated graded module $G(N)$  of $N$ with
respect to the graded left monomial ordering $\prec_{\mathbbm{\sigma
(e)}\hbox{-}gr}$ defined above. \vskip 6pt

{\bf Proof} (i) $\Rightarrow$ (ii) Note that any nonzero homogeneous
element of $G(N)$ is of the form $\sigma (\xi )$ with $\xi \in N$.
If $\G$ is a left Gr\"obner basis of $N$, then there exists some
$g_i\in\G$ such that $\LM (g_i)|\LM (\xi )$, i.e., there is a
monomial $a^{\alpha}\in\B$ such that $\LM (\xi )=\LM (a^{\alpha}\LM
(g_i))$. Since the given left monomial ordering $\PRCEGR$ on $\BE$
is a graded left monomial ordering, it follows from Lemma 2.7 and
Lemma 3.2 that
$$\begin{array}{rcl} \LM (\sigma (\xi )))&=&\sigma (\LM (\xi ))\\
&=&\sigma (\LM (a^{\alpha}\LM (g_i)))\\
&=&\LM (\sigma (a^{\alpha}\LM (g_i)))\\
&=&\LM (\sigma (a)^{\alpha}\sigma (\LM (g_i)))\\
&=&\LM (\sigma (a)^{\alpha}\LM (\sigma (g_i))).\end{array}$$ This
shows that $\LM (\sigma (g_i))|\LM (\sigma (\xi ) )$, thereby
$\sigma (\G )$ is a left Gr\"obner basis for $G(N)$.}\par

(ii) $\Rightarrow$ (i) Suppose that $\sigma (\G )$ is a left
Gr\"obner basis of $G(N)$ with respect to $\prec_{\mathbbm{\sigma
(e)}\hbox{-}gr}$. If $\xi\in N$ and $\xi\ne 0$, then $\sigma (\xi
)\ne 0$, and there exists a $\sigma (g_i)\in\sigma (\G )$ such that
$\LM (\sigma (g_i))|\LM (\sigma (\xi ))$, i.e., there is a monomial
$\sigma (a)^{\alpha}\in\sigma (\B )$ such that $\LM (\sigma (\xi
))=\LM (\sigma (a)^{\alpha}\LM (\sigma (g_i)))$. Again as $\PRCEGR$
is a left graded monomial ordering on $\BE$, by Lemma 2.7 and Lemma
3.2 we have
$$\begin{array}{rcl} \sigma (\LM (\xi ))&=&\LM (\sigma (\xi ))\\
&=&\LM (\sigma (a)^{\alpha}\LM (\sigma (g_i)))\\
&=&\LM (\sigma (a)^{\alpha}\sigma (\LM (g_i)))\\
&=&\LM (\sigma (a^{\alpha}\LM (g_i)))\\
&=&\sigma (\LM (a^{\alpha}\LM (g_i))).\end{array}$$ This shows that
$d_{\rm fil}(\LM (\xi ))=d_{\rm fil}(\LM (a^{\alpha}\LM (g_i)))$.
Since both $\LM (\xi )$ and $\LM (a^{\alpha}\LM (g_i))$ are
monomials in $\BE$, it follows from  the construction of $FL$ and
Lemma 3.1  that $\LM (\xi )=\LM (a^{\alpha}\LM (g_i))$, i.e., $\LM
(g_i)|\LM (\xi )$. This shows that $\G $ is a left Gr\"obner basis
for $N$.\QED\v5

Similarly, in light of Proposition 2.8 we may define an ordering
$\prec_{\widetilde{e}}$ on the $K$-basis $\widetilde{\B (e)}=\{ 
Z^m\T a^{\alpha}\T e_i~|~Z^m\T a^{\alpha}\in\T\B ,~1\le i\le s\}$  
of the $\NZ$-graded free $\T A$-module $\T L=\oplus_{i=1}^s\T A\T 
e_i$ subject to the rule: for $Z^s\T a^{\alpha}\T e_i,Z^t\T 
a^{\beta}\T e_j\in\widetilde{\B (e)}$,
$$Z^s\T a^{\alpha}\T e_i\prec_{\widetilde{e}}Z^t\T a^{\beta}\T e_j
\Longleftrightarrow
a^{\alpha}e_i\prec_{e\hbox{-}gr}a^{\beta}e_j,~\hbox{or}~a^{\alpha}e_i=a^{\beta}e_j~\hbox{and}~s<t.
$$ {\parindent=0pt\par
\def\PRCEGR{\prec_{e\hbox{-}gr}}

{\bf 3.4. Lemma} With the ordering $\prec_{\widetilde{e}}$ defined 
above, the following statements hold.\par

(i)  $\prec_{\widetilde{e}}$  is a  left monomial ordering on 
$\widetilde{\B (e)}$.\par

(ii) (Compare with Corollary 2.5.) $\LM (\T\xi )=\widetilde{\LM (\xi
)}$ holds for all nonzero $\xi\in L$, , where the monomial orderings
used for $\LM (\T\xi )$ and $\LM (\xi )$ are $\prec_{\widetilde{e}}$ 
and $\PRCEGR$ respectively. \vskip 6pt

{\bf Proof} (i) Noticing that the given monomial ordering
$\prec_{gr}$ for $A$ is a graded monomial ordering with respect to a
positive-degree function $d(~)$ on $A$, it follows from  Theorem
2.4(ii) that  $\T A$ is turned into an $\NZ_{\T\B}$-graded solvable
polynomial algebra by using the  monomial ordering $\prec_{_{\T A}}$
defined on $\T\B$ subject to the rule: $\T a^{\alpha}Z^s\prec_{_{\T
A}}\T a^{\beta}Z^t\Longleftrightarrow a^{\alpha}\prec_{gr}a^{\beta},
\hbox{or}~a^{\alpha}=a^{\beta}~\hbox{and}~s<t,\quad
a^{\alpha},a^{\beta}\in\B,$ where the positive-degree function on
$\T A$ is given by  $d(\T{a_i})=d (a_i)$ for $1\le i\le n$, and
$d(Z)=1$. Moreover, since $\T e_i$ is a homogeneous element of
degree $b_i$ in $\T A$, $1\le i\le s$, by Lemma 2.7, it is then
straightforward to verify that $\prec_{\widetilde{e}}$ is a left 
monomial ordering on $\widetilde{B (e)}$.}\par

(ii) Let $\xi =\sum_{i,j}\lambda_{ij}a^{\alpha (i_j)}e_j$, where
$\lambda_{ij}\in K^*$ and $a^{\alpha (i_j)}\in\B$ with $\alpha
(i_j)=(\alpha_{i_{j1}},\ldots ,\alpha_{i_{jn}})\in\NZ^n$. If $d_{\rm
fil}(\xi )=q$, i.e., $\xi\in F_qL-F_{q-1}L$, then by Lemma 3.1, $\LM
(\xi )=a^{\alpha (i_t)}e_t$ for some $t$ such that $d_{\rm
fil}(a^{\alpha (i_t)}e_t)=d(a^{\alpha (i_t)})+b_t=q$. Since
$\prec_{e\hbox{-}gr}$ is a left graded monomial ordering on $\BE$, 
by Lemma 2.7 we have $\T\xi =\lambda_{it}\T a^{\alpha (i_t)}\T 
e_t+\sum_{j\ne t}\lambda_{ij}Z^{q-\ell_{ij}}\T a^{\alpha (i_j)}\T 
e_j$, where $\ell_{ij}=d_{\rm fil}(a^{\alpha (i_j)}e_j)=d(a^{\alpha 
(i_j)})+d_j$. It follows from the definition of  
$\prec_{\widetilde{e}}$  that $\LM (\T\xi )=\T a^{\alpha (i_t)}\T 
e_t=\widetilde{\LM (\xi )}$, as desired.\QED{\parindent=0pt\v5

{\bf 3.5. Theorem} Let $N$ be a submodule of the filtered free
$A$-module $L=\oplus_{i=1}^sAe_i$, where $L$ is equipped with the
filtration $FL=\{ F_qL\}_{q\in\NZ}$ such that $d_{\rm
fil}(e_i)=b_i$, $1\le i\le s$, and let $\prec_{e\hbox{-}gr}$ be a 
graded left monomial ordering on $\BE$. For a subset $\G =\{ 
g_1,\ldots ,g_m\}$ of $N$, the following two statements are 
equivalent.\par

(i) $\G$ is a left Gr\"obner basis of $N$ with respect to
$\PRCEGR$.\par

(ii) Putting $\tau (\G )=\{ \T g_1,\ldots ,\T g_m\}$ and considering
the filtration $FN=\{ F_qN=F_qL\cap N\}_{q\in\NZ}$ of $N$ induced by
$FL$ (see Section 2), $\tau (\G )$ is a left Gr\"obner basis for the
Rees module $\T N$  of $N$ with respect to the  left monomial
ordering $\prec_{\widetilde{e}}$ defined above. \vskip 6pt

{\bf Proof} (i) $\Rightarrow$ (ii) Note that any nonzero homogeneous
element of $\T N$ is of the form $h_q(\xi )$ for some $\xi\in F_qN$
with $d_{\rm fil}(\xi )=q_1\le q$. By Lemma 2.7, $h_q(\xi
)=Z^{q-q_1}\T \xi$. If $\G$ is a left Gr\"obner basis of $N$, then
there exists some $g_i\in\G$ such that $\LM (g_i)|\LM (\xi )$, i.e.,
there is a monomial $a^{\alpha}\in\B$ such that $\LM (\xi )=\LM
(a^{\alpha}\LM (g_i))$. It follows from Lemma 2.7 and Lemma 3.4 that
$$\begin{array}{rcl} \LM (\T\xi )&=&\widetilde{\LM (\xi )}\\
&=&(\LM (a^{\alpha}\LM (g_i)))\widetilde{~}\\
&=&\LM ((a^{\alpha}\LM (g_i))\widetilde{~})\\
&=&\LM (\T a^{\alpha}\widetilde{\LM (g_i)})\\
&=&\LM (\T a^{\alpha}\LM (\T g_i)).\end{array}$$ Hence, noticing the
definition of $\prec_{\widetilde{e}}$ we have
$$\begin{array}{rcl} \LM (h_q(\xi ))&=&\LM (Z^{q-q_1}\T\xi )\\
&=&Z^{q-q_1}\LM (\T\xi )\\
&=&Z^{q-q_1}\LM (\T a^{\alpha}\LM (\T g_i))\\
&=&\LM (z^{q-q_1}\T a^{\alpha}\LM (\T g_i)).\end{array}$$ This shows
that $\LM (\T g_i)|\LM (h_q(\xi ))$, thereby $\tau (\G )$ is a left
Gr\"obner basis of $\T N$.}\par

(ii) $\Rightarrow$ (i)  If $\xi\in N$ and $\xi\ne 0$,  then
$\T\xi\ne 0$ and $\LM (\T\xi )=\widetilde{\LM (\xi )}$ by Lemma 3.4.
Suppose that $\tau (\G )$ is a left Gr\"obner basis of $\T N$ with
respect to  $\prec_{\widetilde{e}}$. Then there exists some $\T 
g_i\in\tau (\G )$ such that  $\LM (\T g_i)|\LM (\T\xi )$, i.e., 
there is a monomial $Z^{m}\T a^{\gamma}\in\widetilde{\B}$ such that 
$\LM (\T\xi )=\LM (Z^{m}\T a^{\gamma}\LM (\T g_i))$. Since the given 
left monomial ordering $\PRCEGR$ on $\BE$ is a graded left monomial 
ordering, it follows from Lemma 2.7, the definition of 
$\prec_{\widetilde{e}}$ and  Lemma 3.2 that
$$\begin{array}{rcl} \widetilde{a^{\alpha}e_j}=\widetilde{\LM (\xi
)}=\LM (\T\xi )&=&\LM (Z^{m}\T a^{\gamma}\LM (\T g_i))\\
&=&Z^{m}(\LM ((a^{\gamma}\LM (g_i))\widetilde{~}))\\
&=&Z^{m}(\LM (a^{\gamma}\LM (g_i)))\widetilde{~}.\end{array}$$
Noticing the discussion on  $\T L$ and the role played by $Z$ given
before Lemma 2.7,   we must have $m=0$, thereby $\LM (\xi )=\LM
(a^{\gamma}\LM (g_i))$. This shows that $\G$ is a left Gr\"obner
basis for $N$.{\parindent=0pt \v5

{\bf Remark} It is known that Gr\"obner bases for ungraded ideals in
both a commutative polynomial algebra and a noncommutative free
algebra can be obtained via computing homogeneous Gr\"obner bases
for graded ideals in the corresponding homogenized (graded) algebras
(cf. [Fr\"ob], [LS], [Li4]). Similarly for a weighted $\NZ$-filtered
solvable polynomial algebra $A$, by using a  (de)homogenization-like
trick  with respect to the central regular element $Z$ in $\T A$,
the discussion on $\T A$ and $\T L$ presented in Section 2 indeed
enables us to obtain left Gr\"obner bases of submodules (left
ideals) in $L$ (in $A$) via computing homogeneous left Gr\"obner
bases of graded submodules (graded left ideals) in $\T L$ (in $\T
A$). Since this topic is beyond the scope of this paper,  we omit
the detailed discussion here. }

\section*{4. F-Bases and Standard Bases with Respect to Good Filtrations}
Let $A=K[a_1,\ldots ,a_n]$ be a weighted $\NZ$-filtered solvable
polynomial algebra with admissible system $(\B ,\prec )$ and the
$\NZ$-filtration  $FA=\{F_pA\}_{p\in\NZ}$ constructed  with respect
to a given positive-degree function $d(~)$ on  $A$ (see Section 2).
In this section, we introduce F-bases and standard bases
respectively for $\NZ$-filtered left  $A$-modules and their
submodules with respect to good filtrations, and we show that any
two minimal F-bases, respectively any two minimal standard bases
have the same number of elements and the same number of elements of
the same filtered degree.  Moreover, we show that a standard basis
for a submodule $N$ of a filtered free $A$-module $L$ can be
obtained via computing a left Gr\"obner basis  of $N$ with respect
to a graded left monomial ordering. All notions, notations and
conventions used before are maintained. \v5

Let $M$ be an $A$-module. Recall that $M$ is said to be an
$\NZ$-filtered $A$-module if $M$ has a filtration $FM=\{
F_qM\}_{q\in\NZ}$, where each $F_qM$ is a $K$-subspace of $M$, such
that $M=\cup_{q\in\NZ}F_qM$, $F_qM\subseteq F_{q+1}M$ for all
$q\in\NZ$, and $F_pAF_qM\subseteq F_{p+q}M$ for all $p,q\in\NZ$.
{\parindent=0pt\v5

{\bf Convention}  Unless otherwise stated, from now on in the
subsequent sections  a filtered $A$-module $M$  is always meant an
$\NZ$-filtered module with a filtration of the type $FM=\{
F_qM\}_{q\in\NZ}$ as described above.}\v5

Let $G(A)$ be the associated graded algebra of $A$, $\T A$ the Rees
algebra of $A$, and $Z$ the homogeneous element of degree 1 in $\T
A_1$ represented by the multiplicative identity 1 of $A$ (see
Section 2). If $M$ is a filtered $A$-module with the filtration
$FM=\{ F_qM\}_{q\in\NZ}$, then  $M$  has the  associated  graded
$G(A)$-module $G(M)=\oplus_{q\in\NZ}G(M)_q$ with $G(M)_0=F_0M$ and
$G(M)_q=F_qM/F_{q-1}M$ for $q\ge 1$, and the Rees module of $M$ is
defined as the graded $\T A$-module $\T M=\oplus_{q\in\NZ}\T M_q$
with each $\T M_q=F_qM$. As with a  filtered free $A$-module in
Section 2, we have  $\T M/Z\T M\cong G(M)$ as graded $G(A)$-modules,
and $\T M/(1-Z)\T M\cong M$ as $A$-modules. Moreover, we may also
define the filtered degree of a nonzero $\xi\in M$, that is, $d_{\rm
fil}(\xi )=q$ if and only if $\xi\in F_qM-F_{q-1}M$. So, actually as
in Section 2, for $\xi\in M$ with $d_{\rm fil}(\xi )=q$, if we write
$\sigma (\xi )$ for the nonzero homogeneous element of degree $q$
represented by $\xi$ in $G(M)_q$, $\T{\xi}$ for the degree-$q$
homogeneous element represented by $\xi$ in $\T M_q$, and
$h_{q'}(\xi )$ for the degree-$q'$ homogeneous element represented
by $\xi$ in $\T M_{q'}$ with $q<q'$, then $d_{\rm fil}(\xi
)=q=d_{\rm gr}(\sigma (\xi ))=d_{\rm gr}(\T\xi )$, and $d_{\rm
gr}(h_{q'}(\xi ))=q'$. \par

With notation as fixed above, the lemma presented below is a version
of ([LVO], Ch.I, Lemma 5.4, Theorem 5.7) for $\NZ$-filtered modules.
{\parindent=0pt\v5

{\bf 4.1. Lemma} Let $M$ be a filtered $A$-module with the
filtration $FM=\{ F_qM\}_{q\in\NZ}$, and $V=\{ v_1,\ldots ,v_m\}$ a
finite subset of nonzero elements in $M$. The following statements
are equivalent:\par

(i)  There is a subset $S=\{ n_1,\ldots ,n_m\}\subset\NZ$ such that
$$F_qM=\sum^m_{i=1}\left (\sum_{p_i+n_i\le q}F_{p_i}A\right )v_i,\quad q\in\NZ ;$$\par

(ii) $G(M)=\sum^m_{i=1}G(A)\sigma (v_i)$;\par

(iii) $\T M=\sum^m_{i=1}\T A\T{v_i}$.\par\QED}\v5

{\parindent=0pt\v5

{\bf 4.2. Definition} Let $M$ be a filtered  $A$-module with the
filtration $FM=\{ F_qM\}_{q\in\NZ}$, and let $V =\{v_1,\ldots,
v_m\}\subset M$ be a finite subset of nonzero elements. If $V$
satisfies one of the equivalent conditions of Lemma 4.1, then we
call $V$ an {\it F-basis} of $M$ with respect to $FM$. }\v5

Let $M$ be a filtered  $A$-module with the filtration $FM=\{
F_qM\}_{q\in\NZ}$. If $V$ is an F-basis of $M$ with respect to $FM$,
then it is necessary to note that{\parindent=1.35truecm\par

\item{(1)} since $M=\cup_{q\in\NZ}F_qM$, it is clear that $V$ is
certainly a generating set of the $A$-module $M$, i.e.,
$M=\sum^m_{i=1}Av_i$;\par

\item{(2)} due to Lemma 4.1(i), the  filtration $FM$ is usually
referred to as a {\it good filtration} of $M$ in the literature
concerning filtered module theory (cf. [LVO]). \par}{\parindent=0pt

Indeed, if an $A$-module $M=\sum^t_{i=1}Au_i$ is finitely generated
by the subset $U =\{u_1,\ldots ,u_t\}$, and if $S=\{n_1,\ldots
,n_t\}$ is an arbitrarily chosen subset of $\NZ$, then $U$ is an
F-basis of $M$ with respect to the good filtration
$FM=\{F_qM\}_{q\in\NZ}$ defined by setting
$$F_qM=\{ 0\}~\hbox{if}~q<\min\{ n_1,\ldots ,n_m\} ;~\hbox{otherwise}~
F_qM=\sum^t_{i=1}\left (\sum_{p_i+n_i\le q}F_{p_i}A\right )u_i,\quad
q\in\NZ .$$  In particular, if  $L=\oplus_{i=1}^sAe_i$ is a filtered
free $A$-module with the filtration $FL=\{ F_qL\}_{q\in\NZ}$ as
constructed in Section 2 such that $d_{\rm fil}(e_i)=b_i$, $1\le
i\le s$, then  $\{ e_1,\ldots ,e_s\}$ is an F-basis of $L$ with
respect to the good filtration $FL$. }{\parindent=0pt\v5

{\bf 4.3. Definition} Let $M$ be a filtered $A$-module with the
filtration $FM=\{ F_qM\}_{q\in\NZ}$, and suppose that $M$ has an
F-basis $V =\{v_1,\ldots, v_m\}$ with respect to $FM$. If any proper
subset of $V$ cannot be an F-basis of $M$ with respect to $FM$, then
we say that $V$ is a {\it minimal F-basis} of $M$ with respect to
$FM$. }\v5

Note that $A$ is a weighted $\NZ$-filtered $K$-algebra such that
$G(A)=\oplus_{p\in\NZ}G(A)_p$ with $G(A)_0=K$, $\T
A=\oplus_{p\in\NZ}\T A_p$ with $\T A_0=K$, while $K$ is a field.  by
Lemma 4.1 and the classical result on graded modules over an
$\NZ$-graded algebra with the degree-0 homogeneous part a field,  we
have immediately the following{\parindent=0pt\v5

{\bf 4.4. Proposition} Let $M$ be a filtered $A$-module with the
filtration $FM=\{ F_qM\}_{q\in\NZ}$, and $V=\{ v_1,\ldots
,v_m\}\subset M$ a subset of nonzero elements.  Then $V$ is  a
minimal F-basis of $M$ with respect to $FM$ if and only if $\sigma
(V)=\{\sigma (v_1),\ldots ,\sigma (v_m)\}$ is a minimal homogeneous
generating set of $G(M)$ if and only if $\tau (V)=\{\T v_1,\ldots
,\T v_m\}$ is a minimal homogeneous generating set of $\T M$. Hence,
any two minimal F-bases of $M$ with respect to  $FM$ have the same
number of elements and the same number of elements of the same
filtered degree.\par\QED}\v5

Let $M$ be an $\NZ$-filtered $A$-module with the filtration $FM=\{
F_qM\}_{q\in\NZ}$, and let $N$ be a submodule of $M$ with the
filtration $FN=\{ F_qN=N\cap F_qM\}_{q\in\NZ}$ induced by $FM$.
Then, as with a filtered free $A$-module in Section 2, the
associated graded $G(A)$-module $G(N)=\oplus_{q\in\NZ}G(N)_q$ of $N$
with $G(N)_q=F_qN/F_{q-1}N$ is a graded submodule of $G(M)$, and the
Rees module $\T N=\oplus_{q\in\NZ}\T N_q$ of $N$ with $\T N_q=F_qN$
is a graded submodule of $\T M$. {\parindent=0pt\v5

{\bf 4.5. Definition}  Let $M$ be a filtered  $A$-module with the
filtration $FM=\{ F_qM\}_{q\in\NZ}$, and let $N$ be a submodule of
$M$. Consider the filtration $FN=\{ F_qN=N\cap F_qM\}_{q\in\NZ}$ of
$N$ induced by $FM$. If $W=\{ \xi_1,\ldots ,\xi_s\}\subset N$ is an
F-basis with respect to $FN$ in the sense of Definition 4.2, then we
call $W$ a {\it standard basis} of $N$. \v5

{\bf Remark} By referring to Lemma 4.1, one may check that our
definition 4.5 of a standard basis coincides with the classical
Macaulay basis provided $A=K[x_1,\ldots ,x_n]$ is the commutative
polynomial $K$-algebra (cf. [KR] Definition 4.2.13, Theorem 4.3.19),
for, taking the $\NZ$-filtration $FA$ with respect to an arbitrarily
chosen positive-degree function $d(~)$ on $A$, there are graded
algebra isomorphisms $G(A)\cong A$ and $\T A\cong K[x_0,x_1,\ldots
,x_n]$, where $d(x_0)=1$ and $x_0$ plays the role that the central
regular element $Z$ of degree 1 does in $\T A$. Moreover, if
two-sided ideals of a weighted $\NZ$-filtered solvable polynomial
algebra  $A$ are considered, then one may see that our definition
4.5 of a standard basis coincides with the standard basis defined in
[Gol].\v5

{\bf 4.6. Definition} Let $M$ be a filtered $A$-module with the
filtration $FM=\{ F_qM\}_{q\in\NZ}$, and  $N$ a submodule of $M$
with the filtration $FN=\{ F_qN=N\cap F_qM\}_{q\in\NZ}$ induced by
$FM$. Suppose that $N$ has a standard basis  $W=\{\xi_1,\ldots
,\xi_m\}$ with respect to $FN$. If any proper subset of $W$ cannot
be a standard basis for $N$ with respect to $FN$, then we call $W$ a
{\it minimal standard basis} of $N$ with respect to $FN$.}\v5

If $N$ is a submodule of a filtered $A$-module $M$ with filtration
$FM$, then since a standard basis of $N$ is defined as an F-basis of
$N$ with respect to the filtration $FN$ induced by $FM$, the next
proposition follows from Proposition 4.4.{\parindent=0pt\v5

{\bf 4.7. Proposition} Let $M$ be a filtered $A$-module with the
filtration $FM=\{ F_qM\}_{q\in\NZ}$, and $N$ a submodule of $M$ with
the induced filtration $FN=\{ F_qN=N\cap F_qM\}_{q\in\NZ}$. A finite
subset of nonzero elements  $W=\{ \xi_1,\ldots ,\xi_s\}\subset N$
is  a minimal standard basis of $N$ with respect to $FN$ if and only
if $\sigma (W)=\{\sigma (\xi_1),\ldots ,\sigma (\xi_m)\}$ is a
minimal homogeneous generating set of $G(N)$ if and only if $\tau
(W)=\{\T \xi_1,\ldots ,\T\xi_m\}$ is a minimal homogeneous
generating set of $\T N$. Hence, any two minimal standard bases of
$N$  have the same number of elements and the same number of
elements of the same filtered degree.\par\QED}\v5

Since $A$, $G(A)$ and $\T A$ are all Noetherian rings, if a filtered
$A$-module $M$ has an F-basis $V$ with respect to a given filtration
$FM$, then  the existence of a standard basis for a submodule $N$ of
$M$ follows immediately from Lemma 4.1. Our next theorem shows that
a standard basis for a submodule $N$ of a filtered free $A$-module
$L$ can be obtained via computing a left Gr\"obner basis  of $N$
with respect to a graded left monomial ordering.{\parindent=0pt\v5
\def\PRCVE{\prec_{\varepsilon\hbox{-}gr}}\def\BV{\B (\varepsilon )}

{\bf 4.8. Theorem} Let  $L=\oplus_{i=1}^sAe_i$ be a filtered free
$A$-module with the filtration $FL=\{ F_qL\}_{q\in\NZ}$ such that
$d_{\rm fil}(e_i)=b_i$, $1\le i\le s$, and let $\prec_{e\hbox{-}gr}$ 
be a graded left monomial ordering on $\BE$ (see Section 3). If $\G 
=\{ g_1,\ldots ,g_m\}\subset L$ is a left Gr\"obner basis for the 
submodule $N=\sum_{i=1}^mAg_i$ of $L$ with respect to $\PRCEGR$, 
then $\G$ is a standard basis for $N$ in the sense of Definition 
4.5.\vskip 6pt
\def\PRCEGR{\prec_{e\hbox{-}gr}}

{\bf Proof} If $\xi\in F_qN=F_qL\cap N$ and $\xi\ne 0$, then $d_{\rm
fil}(\xi )\le q$ and $\xi $ has a left Gr\"obner representation by
$\G$, that is, $\xi =\sum_{i,j}\lambda_{ij}a^{\alpha (i_j)}g_j$,
where $\lambda_{ij}\in K^*$, $a^{\alpha (i_j)}\in\B$ with $\alpha
(i_j)=(\alpha_{i_{j1}},\ldots ,\alpha_{i_{jn}}) \in\NZ^n$,
satisfying $\LM (a^{\alpha (i_j)}g_j)\preceq_{e\hbox{-}gr}\LM (\xi 
)$. Suppose $d_{\rm fil}(g_j)=n_j$, $1\le j\le m$. Since $\PRCEGR$ 
is a graded left monomial ordering on $\BE$, by Lemma 3.1 we may 
assume that $\LM (g_j)=a^{\beta (j)}e_{t_j}$ with $\beta 
(j)=(\beta_{j_1},\ldots ,\beta_{j_n})\in\NZ^n$ and $1\le t_j\le s$, 
such that $d(a^{\beta (j)})+b_{t_j}=n_j$, where $d(~)$ is the given 
positive-degree function on $A$. Furthermore, by the property 
$(\mathbb{P}2)$ presented in Section 1, we have
$$\LM (a^{\alpha (i_j)}g_j)=\LM (a^{\alpha (i_j)}a^{\beta
(j)}e_{t_j})=a^{\alpha (i_j)+\beta (j)}e_{t_j},$$ and it follows
from Lemma 2.2, Lemma 2.7 and Lemma 3.1 that $d(a^{\alpha
(i_j)})+n_j=d(a^{\alpha (i_j)})+d(a^{\beta (j)})+b_{t_j}=d(a^{\alpha
(i_j)+\beta (j)})+b_{t_j}\le q$. Hence $\xi\in \sum_{j=1}^m\left
(\sum_{p_j+n_j\le q}F_{p_j}A\right )g_j$. This shows that
$F_qN=\sum_{j=1}^m\left (\sum_{p_j+n_j\le q}F_{p_j}A\right )g_j$,
i.e., $\G$ is a standard basis for $N$.}\v5

\section*{5. Computation of Minimal F-Bases and Minimal Standard Bases}
Let $A=K[a_1,\ldots ,a_n]$ be a weighted $\NZ$-filtered solvable
polynomial algebra with  admissible system $(\B ,\prec )$ and the
$\NZ$-filtration  $FA=\{F_pA\}_{p\in\NZ}$ constructed  with respect
to a given positive-degree function $d(~)$ on  $A$ (see Section 2).
In this section we show how to algorithmically compute minimal
F-bases for quotient modules of a filtered free left $A$-module $L$,
and how to algorithmically compute minimal standard bases for
submodules of $L$ in the case where a graded left monomial ordering
$\PRCEGR$ on $L$ is employed.  All notions, notations and
conventions used before are maintained.\v5

We start by a little more preparation.  Let $M$ and $M'$ be filtered
$A$-modules with the filtration $FM=\{ F_qM\}_{q\in\NZ}$ and $FM'=\{
F_qM'\}_{q\in\NZ}$ respectively. Recall that an $A$-module
homomorphism $\varphi$: $M\r M'$ is said to be a {\it filtered
homomorphism} if $\varphi (F_qM)\subseteq F_qM'$ for all $q\in\NZ$.
Let $G(A)$ be the associated $\NZ$-graded algebra of $A$ and $\T A$
the Rees algebra of $A$. Then naturally, a filtered homomorphism
$M~\mapright{\varphi}{}~M'$ induces a graded $G(A)$-module
homomorphism $G(M)~\mapright{G(\varphi )}{}~G(M')$, where if $\xi\in
F_qM$ and $\OV\xi =\xi+F_{q-1}M$ is the coset represented by $\xi$
in $G(M)_q=F_qM/F_{q-1}M$, then $G(\varphi )(\OV\xi )=\varphi (\xi
)+F_{q-1}M'\in G(M')_q=F_qM'/F_{q-1}M'$,  and $\varphi$ induces a
graded $\T A$-module homomorphism $\T
M~\mapright{\T\varphi}{}~\T{M'}$, where if $\xi\in F_qM$ and
$h_q(\xi )$ is the homogeneous element of degree $q$ in $\T
M_q=F_qM$, then $\T\varphi (h_q(\xi ))=h_q(\varphi (\xi
))\in\T{M'}_q=F_qM'$. Moreover, if
$M~\mapright{\varphi}{}~M'~\mapright{\psi}{}~M''$ is a sequence of
filtered homomorphisms, then $G(\psi )\circ G(\varphi
)=G(\psi\circ\varphi )$ and $\T\psi\circ\T\varphi
=\widetilde{\psi\circ\varphi}$.
\par

Furthermore, recall that a filtered  homomorphism
$M~\mapright{\varphi}{}~M'$ is called a {\it strict filtered
homomorphism} if $\varphi (F_qM)=\varphi (M)\cap F_qM'$ for all
$q\in\NZ$. Note that if $N$ is a submodule of $M$ and $\OV M=M/N$,
then, considering the induced filtration $FN=\{ F_qN=N\cap
F_qM\}_{q\in\NZ}$ of $N$ and the induced filtration $F(\OV M)=\{
F_q\OV M=(F_qM+N)/N\}_{q\in\NZ}$ of $\OV M$, the inclusion map
$N\hookrightarrow M$ and the canonical map $M\r \OV M$ are strict
filtered  homomorphisms. Concerning strict filtered homomorphisms
and their associated graded homomorphisms, the next  proposition is
quoted from ([LVO], CH.I, Section 4).{\parindent=0pt\v5

{\bf 5.1. Proposition} Given a sequence of filtered homomorphisms
$$N~\mapright{\varphi}{}~M~\mapright{\psi}{}~M',\leqno{(*)}$$
such that $\psi\circ\varphi =0$, the following statements are
equivalent.\par

(i) The sequence $(*)$ is exact and $\varphi$, $\psi$ are strict
filtered homomorphisms.\par

(ii) The sequence $G(N)~\mapright{G(\varphi
)}{}~G(M)~\mapright{G(\psi ) )}{}~G(M')$ is exact.\par

(iii) The sequence $\T N~\mapright{\T\varphi}{}~\T
M~\mapright{\T\psi}{}~\T{M'}$ is exact.\par\QED}\v5

Let $L=\oplus_{i=1}^mAe_i$ be a filtered free $A$-module with the
filtration $FL=\{ F_qL\}_{q\in\NZ}$ such that $d_{\rm
fil}(e_i)=b_i$, $1\le i\le m$. Then as we have noted in Section 4,
$\{ e_1,\ldots ,e_m\}$ is an F-basis of $L$ with respect to the good
filtration $FL$. Let $N$ be a submodule of $L$,  and let the
quotient module $M=L/N$ be equipped with the filtration
$FM=\{F_qM=(F_qL+N)/N\}_{q\in\NZ}$ induced by $FL$. Without loss of
generality, we assume that $\OV e_i\ne 0$ for $1\le i\le m$, where
each $\OV e_i$ is the coset represented by $e_i$ in $M$. Then we see
that $\{ \OV e_1,\ldots ,\OV e_m\}$ is an F-basis of $M$ with
respect to $FM$. {\parindent=0pt\v5

{\bf 5.2. Lemma} Let $M=L/N$ be as fixed above, and let
$N=\sum^s_{j=1}A\xi_j$ be generated by the set of nonzero elements $ 
U =\{\xi_1,\ldots ,\xi_s\}$, where 
$\xi_{\ell}=\sum_{k=1}^sf_{k\ell}e_k$ with $f_{k\ell}\in A$ and 
$d_{\rm fil}( \xi_{\ell})=q_{\ell}$, $1\le \ell\le s$. The following 
statements hold.\par

(i) If for some $j$, $\xi_j$ has a nonzero term $f_{ij}e_i$ such
that $d_{\rm fil}(f_{ij}e_i)=d_{\rm fil}(\xi_j)=q_j$ and the
coefficient $f_{ij}$ is a nonzero constant, say $f_{ij}=1$ without
loss of generality, then for each $\ell =1,\ldots ,j-1,j+1,\ldots
,s$, the element $\xi_{\ell}'=\xi_{\ell}-f_{i\ell}\xi_j$ does not
involve $e_i$. Putting $ U '=\{ \xi_1',\ldots
,\xi'_{j-1},\xi_{j+1}',\ldots ,\xi'_s\}$, $N'=\sum_{\xi_{\ell}'\in U 
'}A\xi_{\ell}'$, and considering the filtered free $A$-module 
$L'=\oplus_{k\ne i}Ae_k$ with the filtration $FL'=\{ 
F_qL'\}_{q\in\NZ}$  in which each $e_k$ has the same filtered degree 
as it is in $L$, i.e., $d_{\rm fil}(e_k)=b_k$, if the quotient 
module $M'=L'/N'$ is equipped with the  filtration $FM'=\{ 
F_qM'=(F_qL'+N')/N'\}_{q\in\NZ}$ induced by $FL'$, then there is a 
strict filtered isomorphism $\varphi$: $M'\cong M$, i.e., $\varphi$ 
is an $A$-module isomorphism such that $\varphi (F_qM')=F_qM$ for 
all $q\in\NZ$. \par

(ii) With the assumptions and notations as in (i), if $ U
=\{\xi_1,\ldots ,\xi_s\}$ is a standard basis of $N$ with respect to
the filtration $FN$ induced by $FL$, then $ U '=\{ \xi_1',\ldots
,\xi'_{j-1},\xi_{j+1}',\ldots ,\xi'_s\}$ is a standard basis of $N'$
with respect to the filtration $FN'$ induced by $FL'$. \vskip 6pt

{\bf Proof} (i) Since $f_{ij}=1$ by the assumption, we see that
every $\xi_{\ell}'=\xi_{\ell}-f_{i\ell}\xi_j$ does not involve
$e_i$. Let $ U '=\{ \xi_1',\ldots ,\xi'_{j-1},\xi_{j+1}',\ldots
,\xi'_s\}$ and  $N'=\sum_{\xi_{\ell}'\in U '}A\xi_{\ell}'$. Then
$N'\subset L'=\oplus_{k\ne i}Ae_k\subset L$ and $N=N'+A\xi_j$. Since
$\xi_j =e_i+\sum_{k\ne i}f_{kj}e_k$, the naturally defined
$A$-module homomorphism $M'=L'/N'\mapright{\varphi}{}L/N=M$ with
$\varphi (\OV e_k)=\OV e_k$, $k=1,\ldots ,i-1,i+1,\ldots ,m$, is an
isomorphism, where, without confusion, $\OV e_k$ denotes the coset
represented by $e_k$ in $M'$ and $M$ respectively. It remains to see
that $\varphi$ is a strict filtered isomorphism. Note that $\{
e_1,\ldots ,e_m\}$ is an F-basis of $L$ with respect to $FL$ such
that $d_{\rm fil}(e_i)=b_i$, $1\le i\le m$, i.e.,
$$F_qL=\sum_{i=1}^m\left (\sum_{p_i+b_i\le q}F_{p_i}A\right )e_i, \quad q\in\NZ ,$$
that $\{ e_1,\ldots ,e_{i-1},e_{i+1},\ldots ,e_m\}$ is an F-basis of
$L'$ with respect to $FL'$ such that $d_{\rm fil}(e_k)=b_k$, where
$k\ne i$, i.e.,
$$F_qL'=\sum_{k\ne i}\left (\sum_{p_i+b_k\le q}F_{p_i}A\right )e_k, \quad q\in\NZ ,$$
and that $\xi_j =e_i+\sum_{k\ne i}f_{kj}e_k$ with $q_j=d_{\rm
fil}(\xi_j)=d_{\rm fil}(e_i)=b_i$ such that $d_{\rm
fil}(f_{kj})+b_k\le q_j$ for all $f_{kj}\ne 0$. It follows that
$\sum_{k\ne i}f_{kj}\OV e_k\in F_{q_j}M'$ and $\varphi (\sum_{k\ne
i}f_{kj}\OV e_k)=\OV e_i\in F_{q_j}M$, thereby $\varphi
(F_qM')=F_qM$ for all $q\in\NZ$, as desired.}\par

(ii) Note that $\xi_{\ell}'=\xi_{\ell}-f_{i\ell}\xi_j$. By the
assumption on $\xi_j$, if $f_{i\ell}\ne 0$ and $d_{\rm
fil}(f_{i\ell}e_i)=d_{\rm fil}(\xi_{\ell})=q_{\ell}$, then since
$d_{\rm fil}(\xi_j)=d_{\rm fil}(e_i)$ we have $d_{\rm
fil}(f_{i\ell}\xi_j)=d_{\rm fil}(\xi_{\ell})=q_{\ell}$. It follows
that if we equip $N$ with the filtration $FN=\{ F_qN=N\cap
F_qL\}_{q\in\NZ}$ induced by $FL$ and consider the associated graded
module $G(N)$ of $N$, then $d_{\rm gr}(\sigma (\xi_{\ell}))=d_{\rm
gr}(\sigma (f_{i\ell}\xi_j))=d_{\rm gr}(\sigma (f_{i\ell})\sigma
(\xi_j))$ in $G(N)$, i.e., $\sigma (\xi_{\ell})-\sigma
(f_{i\ell})\sigma (\xi_j)\in G(N)_{q_{\ell}}$. So, if $\sigma
(\xi_{\ell})-\sigma (f_{i\ell})\sigma (\xi_j)\ne 0$ then $d_{\rm
fil}(\xi_{\ell}')=q_{\ell}$ and thus $$\sigma (\xi_{\ell}')=\sigma
(\xi_{\ell}-f_{i\ell}\xi_j)=\sigma (\xi_{\ell})-\sigma
(f_{i\ell})\sigma (\xi_j).\eqno{(1)}$$ If $f_{i\ell}\ne 0$ and
$d_{\rm fil}(f_{i\ell}e_i)<d_{\rm fil}(\xi_{\ell})=q_{\ell}$, then
$\sigma (\xi_{\ell})=\sum_{d(f_{k\ell})+b_k=q_{\ell}}\sigma
(f_{k\ell})\sigma (e_k)$ does not involve $\sigma (e_i)$. Also since
$d_{\rm fil}(\xi_j)=d_{\rm fil}(e_i)$, we have $d_{\rm
fil}(f_{i\ell}\xi_j)<d_{\rm fil}(\xi_{\ell})=q_{\ell}$. Hence
$d_{\rm fil}(\xi_{\ell}')=d_{\rm fil}(\xi_{\ell})=q_{\ell}$ and thus
$$\sigma (\xi_{\ell}')=\sigma (\xi_{\ell}-f_{i\ell}\xi_j)=\sigma
(\xi_{\ell}).\eqno{(2)}$$ If $f_{i\ell}=0$, then it is clear that we
have the same result as presented in (2). Now, if $ U
=\{\xi_1,\ldots ,\xi_s\}$ is a standard basis of $N$ with respect to
the induced filtration $FN$, then $G(N)=\sum_{\ell =1}^sG(A)\sigma
(\xi_{\ell})$ by Lemma 4.1. But since we have also
$G(N)=\sum_{\xi_{\ell}'\in U '}G(A)\sigma (\xi_{\ell}')+G(A)\sigma 
(\xi_j)$ where the $\sigma (\xi_{\ell}')$ are those nonzero 
homogeneous elements obtained according to the above (1) and (2), it 
follows from Lemma 4.1 that $ U'\cup\{\xi_j\}$ is a standard basis 
of $N$ with respect to the induced filtration $FN$.
\par

We next prove that $ U '$ is a standard basis of
$N'=\sum_{\xi_{\ell}'\in U '}A\xi_{\ell}'$ with respect to the
filtration $FN'=\{ F_qN'=N'\cap F_qL'\}_{q\in\NZ}$ induced by $FL'$.
Since $\{ e_1,\ldots ,e_{i-1},e_{i+1},\ldots ,e_m\}$ is an F-basis
of $L'$ with respect to the filtration $FL'$ such that each $e_k$
has the same filtered degree as it is in $L$, i.e., $d_{\rm
fil}(e_k)=b_k$, it is clear that $F_qL'=L'\cap F_qL$, $q\in\NZ$,
i.e., the filtration $FL'$ is the one induced by $FL$. Considering
the filtration $FN'$ of $N'$ induced by $FL'$, it turns out that
$$F_qN'=N'\cap F_qL'=N'\cap F_qL\subseteq N\cap F_qL=F_qN,\quad q\in\NZ .\eqno{(3)}$$
If $\xi\in F_qN'$, then since $ U '\cup\{\xi_j\}$ is a standard
basis of $N$ with respect to the induced filtration $FN$, the
formula (3) entails that
$$\xi =\sum_{\xi_{\ell}'\in U '}f_{\ell}\xi_{\ell}'+f_j\xi_j~\hbox{with}~f_{\ell},f_j\in A,
~d(f_{\ell})+d_{\rm fil}(\xi_{\ell}')\le q,~d(f_j)+d_{\rm
fil}(\xi_j)\le q.\eqno{(4)}$$ Note that every $\xi_{\ell}'$ does not
involve $e_i$ and consequently $\xi$ does not involve $e_i$. Hence
$f_j=0$ in (4) by the assumption on $\xi_j$, and thus
$\xi\in\sum_{\xi_{\ell}'\in U '}\left (\sum_{p_i+q_{\ell}\le
q}F_{p_i}A\right )\xi_{\ell}'$. Therefor we conclude that $ U '$ is 
a standard basis for $N'$ with respect to the induced filtration 
$FN'$, as desired. \QED \v5

Combining  Proposition 4.4, we now show that an analogue of ([KR],
Proposition 4.7.24) and ([Li4], Proposition 4.2) for quotient
modules of filtered free $A$-modules holds true. {\parindent=0pt\v5

{\bf 5.3. Proposition}  Let $L=\oplus_{i=1}^mAe_i$, $M=L/N$ be as in
Lemma 5.2, and suppose that $ U = \{ \xi_1,...,\xi_s\}$ is now a
standard basis of $N$ with respect to the filtration $FN$ induced by
$FL$. The algorithm presented below computes a subset $\{
e_{i_1},\ldots ,e_{i_{m'}}\}\subset\{ e_1,\ldots ,e_m\}$ and a
subset $V=\{ v_1,\ldots ,v_t\}\subset N\cap L'$, where $m'\le m$ and
$L'=\oplus_{q=1}^{m'}Ae_{i_q}$ such that\par

(i) there is a strict filtered isomorphism $L'/N'=M'\cong M$, where
$N'=\sum^{t}_{\ell =1}Av_{\ell}$, and $M'$ has the filtration
$FM'=\{ F_qM'=(F_qL'+N')/N'\}_{q\in\NZ}$ induced by the good
filtration $FL'$ determined by the F-basis $\{ e_{i_1},\ldots
,e_{i_{m'}}\}$ of $L'$;\par

(ii) $V=\{ v_1,\ldots ,v_t\}$ is a standard basis of
$N'=\sum^{t}_{\ell =1}Av_{\ell}$ with respect to the filtration
$FN'$ induced by $FL'$, such that each
$v_{\ell}=\sum^{m'}_{k=1}h_{k\ell}e_{i_k}$ has the property that
$h_{k\ell}\in K^*$ implies $d_{\rm fil}(e_{i_k})=b_{i_k}<d_{\rm
fil}(v_{\ell})$;\par

(iii) $\{\OV e_{i_1},\ldots ,\OV e_{i_{m'}}\}$ is a minimal F-basis
of $M$ with respect to the filtration $FM$. \vskip 6pt

\underline{\bf Algorithm 2
~~~~~~~~~~~~~~~~~~~~~~~~~~~~~~~~~~~~~~~~~~~~~~~~~~~~~~~~~~~~~~~~~~~~~~~~~~~~~~~~~~~~~~~~~~~~~~~~~}\par

$\begin{array}{l} \textsc{INPUT}: E=\{ e_1,\ldots ,e_m\};~~ U = \{ 
\xi_1,...,\xi_s\},~ 
\hbox{where}~\xi_{\ell}=\sum_{k=1}^mf_{k\ell}e_k~
\hbox{with}~f_{k\ell}\in A,\\
~~~~~~~~~~~~~~~~~~~~~~~~~~~~~~~~~~~~~~~~~~~~~~~~~~~~~~~~~~~~\hbox{and}~d(f_{k\ell})+b_k\le q_{\ell}=d_{\rm fil}(\xi_{\ell}), ~1\le \ell\le s.\\
\textsc{OUTPUT}: ~E' =\{ e_{i_1},\ldots ,e_{i_{m'}}\}; ~~V=\{
v_1,\ldots ,v_t\}\subset N\cap L',~\hbox{where}~
L'=\oplus_{k=1}^{m'}Ae_{i_k}\\
\textsc{INITIALIZATION}:    ~t :=s;~ V := U ;  ~m':=m; ~E' :=E;
\end{array}$\par

$\begin{array}{l}
\textsc{BEGIN}\\
~~~~\textsc{WHILE}~\hbox{there is a}~v_j=\sum^{m'}_{k=1}f_{kj}e_k\in
V~\hbox{and}~i~
\hbox{is the least index}\\
~~~~~~~~~~~~~~~~~\hbox{such that}~f_{ij}\in K^*~
\hbox{with}~d(f_{ij})+b_i=d_{\rm fil}(v_j)~\hbox{DO}\\
~~~~~~~~~~~~~~~~~\hbox{set}~T=\{1,\ldots ,j-1,j+1,\ldots
,t\}~\hbox{and compute}\\
~~~~~~~~~~~~~~~~~v_{\ell}'
=v_{\ell}-\frac{1}{f_{ij}}f_{i\ell}v_j,~\ell\in T,\\
~~~~~~~~~~~~~~~~~r=\#\{\ell~|~\ell\in T,~v_{\ell} =0\}\\
~~~~~~~~~~t := t-r-1\\
~~~~~~~~~V :=\{ v_{\ell}=v_{\ell}'~|~\ell\in T,~v_{\ell}'\ne 0\}\\
~~~~~~~~~~~~=\{ v_1,\ldots ,v_t\}~(\hbox{after reordered})\\
~~~~~~~~~m' :=m'-1\\
~~~~~~~~~E':=E'-\{ e_i\}=\{ e_1,\ldots ,e_{m'}\}~(\hbox{after reordered})\\
~~~~\textsc{END}\\
\textsc{END}\end{array}$\par
\underline{~~~~~~~~~~~~~~~~~~~~~~~~~~~~~~~~~~~~~~~~~~~~~~~~~~~~~~~~~~~~~~~~~~~~~~~~~~~~~~~~~~~~~~~~~~~~~~~~~~~~~~~~~~~~~~~~~~~~~~~~~~~~~~~~~~}
\vskip 6pt

{\bf Proof} First note that for each $\xi_{\ell}\in U$, $d_{\rm
fil}(\xi_{\ell})$ is determined by Lemma 2.6. Since the algorithm is
clearly finite, the conclusions (i) and (ii) follow from Lemma
5.2.}\par

To prove the conclusion (iii), by the strict filtered isomorphism
$M'=L'/N'\cong M$ (or the proof of Lemma 5.2(i)) it is sufficient to
show that $\{\OV e_{i_1},\ldots ,\OV e_{i_{m'}}\}$ is a minimal
F-basis of $M'$ with respect to the filtration $FM'$. By the
conclusion (ii), $V=\{ v_1,\ldots ,v_t\}$ is a standard basis of the
submodule $N'=\sum^{t}_{\ell =1}Av_{\ell}$ of $L'$ with respect to
the filtration $FN'$ induced by $FL'$ such that each
$v_{\ell}=\sum^{m'}_{k=1}h_{k\ell}e_{i_k}$ has the property that
$h_{k\ell}\in K^*$ implies $d_{\rm fil}(e_{i_k})=b_{i_k}<d_{\rm
fil}(v_{\ell})$. It follows from Lemma 4.1 that
$G(N')=\sum^t_{k=1}G(A)\sigma (v_k)$ in which each $\sigma
(v_k)=\sum_{d(h_{k\ell})+b_{i_k}=d_{\rm fil}(v_k)}\sigma
(h_{k\ell})\sigma (e_{i_k})$ and all the coefficients $\sigma
(h_{k\ell})$ satisfy $d_{\rm gr}(\sigma (h_{k\ell}))>0$ (see Section
2). Noticing that $G(A)_0=K$, $G(L')=\oplus_{k=1}^{m'}G(A)\sigma
(e_{i_k})$ (Proposition 2.8) and $G(N')$ is the graded syzygy module
of the graded quotient module $G(L')/G(N')$, by the classical result
on finitely generated graded modules over $\NZ$-graded algebras with
the degree-0 homogeneous part a field,  we conclude  that $\{
\OV{\sigma (e_{i_1})},\ldots ,\OV{\sigma (e_{i_{m'}})}\}$ is a
minimal homogeneous generating set of $G(L')/G(N')$. On the other
hand, considering the naturally formed exact sequence of  strict
filtered homomorphisms
$$0~\mapright{}{}~N'~\mapright{}{}~L'~\mapright{}{}~M'=L'/N'~\mapright{}{}~0,$$
by Proposition 5.1 we have the $\NZ$-graded $G(A)$-module
isomorphism $G(L')/G(N')\cong G(L'/N')=G(M')$ with $\OV{\sigma
(e_{i_k})}\mapsto \sigma (\OV e_{i_k})$, $1\le k\le m'$. Now,
applying Proposition 4.4, we conclude that $\{\OV e_{i_1},\ldots
,\OV e_{i_{m'}}\}$ is a minimal F-basis of $M'$ with respect to the
filtration $FM'$, as desired. \QED \v5

Finally, let $L=\oplus_{i=1}^sAe_i$ be a filtered free $A$-module
with the filtration $FL=\{ F_qL\}_{q\in\NZ}$ such that $d_{\rm
fil}(e_i)=b_i$, $1\le i\le s$, and let $\prec_{e\hbox{-}gr}$ be a 
graded left monomial ordering on the $K$-basis $\BE =\{ 
a^{\alpha}e_i~|~a^{\alpha}\in\B ,~1\le i\le s\}$ of $L$ (see Section 
3). Combining Theorem 3.3 and ([Li3], Theorem 3.8), the next theorem 
shows how to algorithmically compute a minimal standard basis. 
{\parindent=0pt\v5

{\bf 5.4. Theorem} Let $N=\sum^c_{i=1}A\theta_i$ be a submodule of
$L$ generated by the subset of nonzero elements $\Theta =\{
\theta_1,\ldots ,\theta_c\}$, and let $FN=\{ F_qN=F_qL\cap
N\}_{q\in\NZ}$  be the filtration of $N$ induced by $FL$. Then a
minimal standard basis of  $N$ with respect to $FN$ can be obtained
by implementing the following procedures:}\par

{\bf Procudure 1.} Run {\bf Algorithm 1} presented in Section 1 with
the initial input data $\Theta =\{ \theta_1,\ldots ,\theta_c\}$ to
compute a left Gr\"obner basis $ U =\{ \xi_1,\ldots ,\xi_m\}$ for
$N$ with respect to $\prec_{e\hbox{-}gr}$ on $\BE$.
\par

{\bf Procudure 2.} Let $G(N)$ be the associated graded $G(A)$-module
of $N$ determined by the induced filtration $FN$. Then $G(N)$ is a
graded submodule of the associated grade free $G(A)$-module $G(L)$,
and it  follows from  Theorem 3.3 that $\sigma ( U )=\{\sigma
(\xi_1),\ldots ,\sigma (\xi_m)\}$ is a homogeneous left Gr\"obner
basis of $G(N)$ with respect to $\prec_{\sigma (e)\hbox{-}gr}$ on 
$\sigma (\BE )$. Run {\bf Algorithm 3} presented in ([Li3], Theorem 
3.8)  with the initial input data $\sigma ( U )$ to compute a 
minimal homogeneous generating set $\{ \sigma (\xi_{j_1}),\ldots 
,\sigma (\xi_{j_t})\}\subseteq\sigma ( U )$ for $G(N)$.
\par

{\bf Procudure 3.} Write down $W=\{ \xi_{j_1},\ldots ,\xi_{j_t}\}$.
Then $W$ is a Minimal standard basis for $N$ by Proposition 4.7.
\par\QED{\parindent=0pt\v5

{\bf Remark} By Theorem 3.5 and Proposition 4.7 it is clear that we
can also obtain a minimal standard basis of the submodule $N$ via
computing a minimal homogeneous generating set for the Rees module
$\T N$ of $N$, which is a graded submodule of the Rees module $\T L$
of $L$. However, noticing the structure of $\T A$ and $\T L$ (see
Theorem 2.4, Proposition 2.8) it is equally clear that the cost of
working on $\T A$ will be much higher than working on $G(N)$.}\v5

\section*{6. The Uniqueness of Minimal Filtered Free Resolutions}
Let $A$ be a weighted $\NZ$-filtered solvable polynomial algebra
with admissible system $(\B ,\prec )$ and the $\NZ$-filtration
$FA=\{F_pA\}_{p\in\NZ}$ constructed  with respect to a given
positive-degree function $d(~)$ on  $A$ (see Section 2).  In this
section, by using minimal F-bases  and minimal standard bases in the
sense of Definition 4.3 and Definition 4.6, we define minimal
filtered free resolutions for finitely generated left $A$-modules,
and we show that such minimal resolutions are  unique up to strict
filtered isomorphism of chain complexes in the category of filtered
$A$-modules. All notions, notations and conventions used before are
maintained.\v5

Let $M=\sum^m_{i=1}A\xi_i$ be an arbitrary finitely generated
$A$-module. Then,  as we have noted in Section 4,  $M$ may be
endowed with a good filtration $FM=\{ F_qM\}_{q\in\NZ}$ with respect
to an arbitrarily chosen subset $\{ n_1,\ldots ,n_m\}\subset\NZ$,
where
$$F_qM=\{ 0\}~\hbox{if}~q<\min\{ n_1,\ldots ,n_m\} ;~\hbox{otherwise}~
F_qM=\sum^t_{i=1}\left (\sum_{p_i+n_i\le q}F_{p_i}A\right
)\xi_i,\quad q\in\NZ .$$ If we consider the filtered free $A$-module
$L_0=\oplus_{i=1}^mAe_i$ with the good filtration $FL_0=\{
F_qL_0\}_{q\in\NZ}$ such that $d_{\rm fil}(e_i)=n_i$, $1\le i\le m$,
then it follows from the construction of $FL_0$ (see Section 2) and
Proposition 5.1 that the following proposition
holds.{\parindent=0pt\v5

{\bf 6.1. Proposition} (i) There is an exact sequence of strict
filtered homomorphisms
$$0~\mapright{}{}~N~\mapright{\iota}{}~L_0~\mapright{\varphi_0}{}~M~\mapright{}{}~0,$$
in which $\varphi_0(e_i)=\xi_i$, $1\le i\le m$, $N=$ Ker$\varphi_0$
with the induced filtration $FN=\{ F_qN=N\cap F_qL_0\}_{q\in\NZ}$,
and $\iota$ is the inclusion map. Equipping $\OV L_0=L_0/N$ with the
induced filtration $F\OV L_0=\{ F_q\OV
L_0=(F_qL_0+N)/N\}_{q\in\NZ}$, it turns out that the induced
$A$-module isomorphism $\OV L_0~\mapright{\OV\varphi_0}{\cong}~M$ is
a strict filtered  isomorphism, that is, $\OV\varphi_0$ satisfies
$\OV\varphi_0(F_q\OV L_0)=F_qM$ for all $q\in\NZ$.
\par

(ii) The induced sequence
$$0~\mapright{}{}~G(N)~\mapright{G(\iota )}{}~G(L_0)~\mapright{G(\varphi_0)}{}~G(M)~\mapright{}{}~0$$
is an exact sequence of graded $G(A)$-modules homomorphisms, thereby
$G(L_0)/G(N)\cong G(M)\cong G(L_0/N)$ as graded $G(A)$-modules. \par

(iii) The induced sequence
$$0~\mapright{}{}~\T N~\mapright{\T\iota}{}~\T{L_0}~\mapright{\T{\varphi_0}}{}~\T M~\mapright{}{}~0$$
is an exact sequence of graded $\T A$-modules  homomorphisms,
thereby $\T{L_0}/\T N\cong \T M\cong\widetilde{L_0/N}$ as graded $\T
A$-modules.\par\QED}\v5

Proposition 6.1(i) enables us to make the
following{\parindent=0pt\v5

{\bf Convention} In what follows we shall always assume that a
finitely generated $A$-module $M$ is of the form as presented in
Proposition 5.1(i), i.e., $M=L_0/N$ and $M$ has the good filtration
$FM=\{F_qM=(F_qL_0+N)/N\}_{q\in\NZ}$.}\v5

Comparing with the classical minimal graded free resolutions defined
for finitely generated graded modules over $\NZ$-graded  algebras
with the degree-0 homogeneous part a field, the results obtained in
previous sections and the preliminary we made above naturally
motivate the following{\parindent=0pt\v5

{\bf 6.2. Definition} Let $L_0=\oplus_{i=1}^mAe_i$ be a filtered
free $A$-module with the filtration $FL_0=\{ F_qL_0\}_{q\in\NZ}$
such that $d_{\rm fil}(e_i)=b_i$, $1\le i\le m$, let $N$ be a
submodule of $L_0$,  and let the $A$-module $M=L_0/N$ be equipped
with the filtration $FM=\{F_qM=(F_qL_0+N)/N\}_{q\in\NZ}$. A {\it
minimal filtered free resolution} of $M$ is an exact sequence of
filtered $A$-module  homomorphisms
$${\cal L}_{\bullet}\quad\quad\cdots~\mapright{\varphi_{i+1}}{}~L_i~\mapright{\varphi_i}{}~\cdots ~
\mapright{\varphi_2}{}~L_1~\mapright{\varphi_1}{}~L_0~\mapright{\varphi_0}{}~M~\r~0$$
satisfying{\parindent=1.3truecm\par

\item{(1)} $\varphi_0$ is the canonical epimorphism, i.e.,
$\varphi_0(e_i)=\OV e_i$ for $e_i\in \mathscr{E}_0=\{ e_1,\ldots
,e_m\}$ (where each $\OV e_i$ denotes the coset represented by $e_i$
in $M$), such that $\varphi_0(\mathscr{E}_0)$ is a minimal F-basis
of $M$ with respect to  $FM$ (in the sense of Definition 4.3);
\par

\item{(2)} for $i\ge 1$,  each $L_i$ is a filtered free
$A$-module with finite $A$-basis $\mathscr{E}_i$,  and each
$\varphi_i$ is a strict filtered homomorphism, such that
$\varphi_i(\mathscr{E}_i)$ is a minimal standard basis of
Ker$\varphi_{i-1}$ with respect to the filtration induced by
$FL_{i-1}$ (in the sense of Definition 4.6). \v5}}

To see that the minimal filtered free resolution introduced above is
an appropriate definition for finitely generated modules over
weighted $\NZ$-filtered solvable polynomial algebras, we  now show
that such a resolution is unique up to a strict filtered isomorphism
of chain complexes in the category of filtered $A$-modules.
{\parindent=0pt\v5

{\bf 6.3. Theorem} Let ${\cal L}_{\bullet}$ be a minimal filtered
free resolution of $M$ as presented in Definition 6.2. The following
statements hold.\par

(i) The associated sequence of graded $G(A)$-modules and graded
$G(A)$-module homomorphisms
$$G({\cal L}_{\bullet})\quad\quad\cdots~\mapright{G(\varphi_{i+1})}{}~G(L_i)~\mapright{G(\varphi_i)}{}~\cdots ~
\mapright{G(\varphi_2)}{}~G(L_1)~\mapright{G(\varphi_1)}{}~G(L_0)~\mapright{G(\varphi_0)}{}~G(M)~\r~0$$
is a minimal graded free resolution of the finitely generated graded
$G(A)$-module $G(M)$.\par

(ii) The associated sequence of graded $\T A$-modules and graded $\T
A$-module homomorphisms
$$\widetilde{\cal L}_{\bullet}\quad\quad\cdots~\mapright{\T\varphi_{i+1}}{}~\T L_i~\mapright{\T\varphi_i}{}
~\cdots ~ \mapright{\T\varphi_2}{}~\T
L_1~\mapright{\T\varphi_1}{}~\T L_0~\mapright{\T\varphi_0}{}~\T
M~\r~0$$ is a minimal graded free resolution of the finitely
generated graded $\T A$-module $\T M$.

(iii) ${\cal L}_{\bullet}$ is uniquely determined by $M$ in the
sense that if $M$ has another minimal filtered free resolution
$${\cal L'}_{\bullet}\quad\quad\cdots~\mapright{\varphi'_{i+1}}{}~L'_i~\mapright{\varphi'_i}{}
~\cdots ~
\mapright{\varphi'_2}{}~L_1'~\mapright{\varphi'_1}{}~L_0~\mapright{\varphi_0}{}~M~\r~0$$
then for each $i\ge 1$, there is a strict filtered $A$-module
isomorphism  $\chi_i$: $L_i\r L_i'$ such that the diagram
$$\begin{array}{ccc} L_i&\mapright{\varphi_i}{}&L_{i-1}\\
\mapdown{\chi_i}\scriptstyle{\cong}&&\mapdown{\chi_{i-1}}\scriptstyle{\cong}\\
L_i'&\mapright{\varphi_i'}{}&L_{i-1}'\end{array}$$ is commutative,
thereby $\{\chi_i~|~i\ge 1\}$ gives rise to a strict filtered
isomorphism of chain complexes of filtered modules: ${\cal
L}_{\bullet}\cong {\cal L}_{\bullet}'$. \vskip 6pt

{\bf Proof} (i) and (ii) follow from Proposition 4.4, Proposition
4.7, and Proposition 6.1.}\par

To prove (iii), let the sequence ${\cal L}_{\bullet}'$ be as
presented above. By (ii), $G({\cal L}_{\bullet}')$ is another
minimal graded free resolution of $G(M)$. It follows from the
classical result on minimal graded free resolutions that there is a
graded isomorphism of chain complexes $G({\cal L}_{\bullet})\cong
G({\cal L}_{\bullet}')$ in the category of graded $G(A)$-modules,
i.e., for each $i\ge 1$, there is a graded $G(A)$-modules
isomorphism  $\psi_i$: $G(L_i)\r G(L_i')$ such that the diagram
$$\begin{array}{ccc} G(L_i)&\mapright{G(\varphi_i)}{}&G(L_{i-1})\\
\mapdown{\psi_i}\scriptstyle{\cong}&&\mapdown{\psi_{i-1}}\scriptstyle{\cong}\\
G(L_i')&\mapright{G(\varphi_i')}{}&G(L_{i-1}')\end{array}$$ is
commutative. Our aim below is to construct the desired strict
filtered isomorphisms $\chi_i$ by using the graded isomorphisms
$\psi_i$ carefully. So, starting with $L_0$, we assume that we have
constructed the strict filtered isomorphisms
$L_j~\mapright{\chi_j}{}~L_j'$, such that $G(\chi_j)=\psi_j$ and
$\chi_{j-1}\varphi_j=\varphi_j'\chi_j$, $1\le j\le i-1$. Let
$L_i=\oplus_{j=1}^{s_i}Ae_{i_j}$. Since each $\psi_i$ is a graded
isomorphism, we have $\psi_i(\sigma (e_{i_j}))=\sigma (\xi_j')$ for
some $\xi_j'\in L_i'$ satisfying $d_{\rm fil}(\xi_j')=d_{\rm
gr}(\sigma (\xi_j'))=d_{\rm gr}(\sigma (e_{i_j}))=d_{\rm
fil}(e_{i_j})$. It follows that if we construct the  filtered
$A$-module homomorphism $L_i~\mapright{\tau_i}{}~L_i'$ by setting
$\tau_i(e_{i_j})=\xi_j'$, $1\le j\le s_i$, then $G(\tau_i)=\psi_i$.
Hence, $\tau_i$ is a strict filtered isomorphism by Proposition 5.1.
Since $\psi_{i-1}=G(\chi_{i-1})$, $\psi_i=G(\tau_i)$, and thus
$$\begin{array}{l}
\psi_{i-1}G(\varphi_i)=G(\chi_{i-1})G(\varphi_i)=G(\chi_{i-1}\varphi_i)\\
G(\varphi_i')\psi_i=G(\varphi_i')G(\tau_i)=G(\varphi_i'\tau_i),\end{array}$$
for each $q\in\NZ$, by the strictness of the $\varphi_j$,
$\varphi_j'$, $\chi_j$ and $\tau_i$, we have
$$\begin{array}{rcl} G(\chi_{i-1}\varphi_i)(G(L_i)_q)&=&((\chi_{i-1}\varphi_i)(F_qL_i)+F_{q-1}L_{i-1}')/F_{q-1}L_{i-1}'\\
&=&(\chi_{i-1}(\varphi_i(L_i)\cap
F_qL_{i-1})+F_{q-1}L_{i-1}')/F_{q-1}L_{i-1}'\\
&\subseteq&((\chi_{i-1}\varphi_i)(L-i)\cap\chi_{i-1}(F_qL_{i-1})+F_{q-1}L_{i-1}')/F_{q-1}L_{i-1}'\\
&=&((\chi_{i-1}\varphi_i)(L_i)\cap
F_qL_{i-1}'+F_{q-1}L_{i-1}')/F_{q-1}L_{i-1}',\end{array}$$

$$\begin{array}{rcl} G(\varphi_i'\tau_i)(G(L_i)_q)&=&((\varphi_i'\tau_i)(F_qL_i)+F_{q-1}L_{i-1}')/F_{q-1}L_{i-1}'\\
&=&(\varphi_i'(F_qL_i')+F_{q-1}L_{i-1}')/F_{q-1}L_{i-1}'\\
&=&(\varphi_i'(L_i')\cap
F_qL_{i-1}'+F_{q-1}L_{i-1}')/F_{q-1}L_{i-1}'\\
&=&((\varphi_i'\tau_i)(L_i)\cap
F_qL_{i-1}'+F_{q-1}L_{i-1}')/F_{q-1}L_{i-1}'.\end{array}$$ Note that
by the exactness of ${\cal L}_{\bullet}$ and ${\cal L}_{\bullet}'$,
as well as the commutativity
$\chi_{i-2}\varphi_{i-1}=\varphi_{i-1}'\chi_{i-1}$, we have
$(\chi_{i-1}\varphi_i)(L_i)\subseteq
\varphi_i'(L_i')=(\varphi_i'\tau_i)(L_i)$. Considering the
filtration of the submodules $(\chi_{i-1}\varphi_i)(L_i)$ and
$(\varphi_i'\tau_i)(L_i)$ induced by the filtration $FL_{i-1}'$ of
$L_{i-1}'$, the commutativity
$\psi_{i-1}G(\varphi_i)=G(\varphi_i')\psi_i$ and the formulas
derived above show that both submodules have the same associated
graded module, i.e.,
$G((\chi_{i-1}\varphi_i)(L_i))=G((\varphi_i'\tau_i)(L_i))$. It
follows from a similar proof of  ([LVO], Theorem 5.7 on P.49) that
$$(\chi_{i-1}\varphi_i)(L_i)=(\varphi_i'\tau_i)(L_i).\eqno{(1)}$$ Clearly, the equality $(1)$
does not necessarily mean the commutativity of the diagram
$$\begin{array}{ccc} L_i&\mapright{\varphi_i}{}&L_{i-1}\\
\mapdown{\tau_i}\scriptstyle{\cong}&&\mapdown{\chi_{i-1}}\scriptstyle{\cong}\\
L_i'&\mapright{\varphi_i'}{}&L_{i-1}'\end{array}$$ To remedy this
problem, we need to further modify the filtered isomorphism
$\tau_i$. Since $G(\chi_{i-1}\varphi_i)(\sigma
(e_{i_j}))=G(\varphi_i'\tau_i)(\sigma (e_{i_j}))$, $1\le j\le s_i$,
if $d_{\rm fil}(e_{i_j})=b_j$, then by (1) and the strictness of
$\tau_i$ and $\varphi_i$ we have
$$\begin{array}{rcl} (\chi_{i-1}\varphi_i)(e_{i_j})-(\varphi_i'\tau_i)(e_{i_j})&\in&
(\varphi_i'\tau_i)(L_i)\cap F_{b_j-1}L_{i-1}'\\
&=&\varphi_i'(L_i')\cap F_{b_j-1}L_{i-1}'\\
&=&\varphi_i'(F_{b_j-1}L_i')\\
&=&(\varphi_i'\tau_i)(F_{b_j-1}L_i),\end{array}\eqno{(2)}$$ and
furthermore from (2) we have a $\xi_j\in F_{b_j-1}L_i$ such that
$d_{\rm fil}(e_{i_j}-\xi_j)=b_j$ and
$$(\chi_{i-1}\varphi_i)(e_{i_j})=(\varphi_i'\tau_i)(e_{i_j}-\xi_j), ~~~
1\le j\le s_i.\eqno{(3)}$$ Now, if we construct the filtered
homomorphism $L_i~\mapright{\chi_i}{}~L_i'$ by setting
$\chi_i(e_{i_j})=\tau_i(e_{i_j}-\xi_j)$, $1\le j\le s_i$, then since
$\tau_i(\xi_j)\in F_{b_j-1}L_i'$, it turns out that
$$G(\chi_i)(\sigma (e_{i_j}))=G(\tau_i)(\sigma (e_{i_j}-\xi_j))=G(\tau_i)(\sigma (e_{i_j}))=\psi_i(\sigma (e_{i_j})),
~~~1\le j\le s_i, $$ thereby $G(\chi_i)=\psi_i$. Hence, $\chi_i$ is
a strict filtered isomorphism by Proposition 5.1, and moreover, it
follows from (3) that we have reached the following diagram
$$\begin{array}{ccccc} \cdots~\mapright{\varphi_{i+1}}{}&L_i&\mapright{\varphi_i}{}
&L_{i-1}&\mapright{\varphi_{i-1}}{}~\cdots\\
&\mapdown{\chi_i}\scriptstyle{\cong}&&\mapdown{\chi_{i-1}}\scriptstyle{\cong}&\\
\cdots~\mapright{\varphi_{i+1}'}{}&L_i'&\mapright{\varphi_i'}{}&L_{i-1}'&\mapright{\varphi_{i-1}'}{}~\cdots\end{array}$$
in which $\chi_{i-1}\varphi_i=\varphi_i'\chi_i$. Repeating the same
process to getting the desired $\chi_{i+1}$ and so on, the proof is
thus finished. \v5

\section*{7. Computation of  Minimal Finite Filtered
Free Resolutions} Let  $A=K[a_1,\ldots ,a_n]$ be a solvable
polynomial algebra with the  admissible system $(\B ,\prec_{gr})$ in
which $\prec_{gr}$ is a graded monomial ordering with respect to
some  given positive-degree function $d(~)$ on  $A$ (see Section
2). Thereby $A$ is turned into a weighted $\NZ$-filtered solvable
polynomial algebra with  the filtration $FA=\{ F_pA\}_{p\in\NZ}$
constructed with respect to the same $d(~)$ (see Example (2) of
Section 2). Our aim in this section is to give algorithmic
procedures for computing  minimal finite filtered free resolutions
over $A$ (in the sense of Definition 6.2). All notions, notations
and conventions used before are maintained.\v5

We start by an arbitrary free left $A$-module $L=\oplus_{i=1}^mAe_i$
with  a left monomial ordering $\prec_{e}$ on the $K$-basis $\BE$ of 
$L$. As in Section 1 we write $S_{\ell}(\xi_i,\xi_j)$ for the left 
S-polynomial of two elements $\xi_i$, $\xi_j\in L$. Let 
$N=\sum^m_{i=1}A\xi_i$ be a submodule of $L$ generated by the set of 
nonzero elements $ U =\{\xi_1,\ldots ,\xi_m\}$. We first demonstrate 
how to calculate a generating set of the syzygy module Syz$( U )$  
by means of a left Gr\"obner basis of $N$. To this end, let $\G=\{ 
g_1,\ldots ,g_t\}$ be a left Gr\"obner basis of $N$ produced by 
running {\bf Algorithm 1} (presented in Section 1) with the initial 
input data $ U$ and the ordering  $\prec_{e}$, then every nonzero 
left S-polynomial $S_{\ell}(g_i,g_j)$ has a left Gr\"obner 
representation $S_{\ell}(g_i,g_j)=\sum_{i=1}^tf_ig_i$ with 
$\LM(f_ig_i)\preceq_{e}\LM (S_{\ell}(g_i,g_j))$ whenever $f_i\ne 0$ 
(note that such a representation is obtained by using the division 
by $\G$ during executing the WHILE loop in {\bf Algorithm 1}). 
Considering the syzygy module Szy$(\G )$ of $\G$ in the free 
$A$-module $L_1=\oplus_{i=1}^tA\varepsilon_i$, if we put
$$s_{ij}=f_1\varepsilon_1+\cdots +\left (f_i-\frac{a^{\gamma-\alpha (i)}}{\LC (a^{\gamma -\alpha (i)}\xi_i)}\right )\varepsilon_i+
\cdots +\left (f_j+\frac{a^{\gamma -\alpha (j)}}{\LC (a^{\gamma
-\alpha (j)}\xi_j)}\right )\varepsilon_j+\cdots +f_t\varepsilon_t,$$
${\cal S} =\{ s_{ij}~|~1\le i<j\le t\}$, then it can be shown,
actually as in the commutative case (cf. [AL], Theorem 3.7.3), that
${\cal S}$ generates Szy$(\G )$. However, by employing an analogue 
of the Schreyer ordering $\prec_{s\hbox{-}\varepsilon}$ on the 
$K$-basis $\B (\varepsilon )=\{ 
a^{\alpha}\varepsilon_i~|~a^{\alpha}\in\B ,~1\le i\le m\}$ of $L_1$ 
induced by $\G$ with respect to $\prec_{e}$, which is defined 
subject to the rule: for 
$a^{\alpha}\varepsilon_i,a^{\beta}\varepsilon_j\in\B (\varepsilon 
)$,
$$a^{\alpha}\varepsilon_i\prec_{s\hbox{-}\varepsilon} a^{\beta}\varepsilon_j
\Leftrightarrow\left\{\begin{array}{l}
\LM (a^{\alpha}g_i)\prec_{e}\LM (a^{\beta}g_j),\\
\hbox{or}\\
\LM (a^{\alpha}g_i)=\LM (a^{\beta}g_j)~\hbox{and}~i<j,
\end{array}\right.$$
there is indeed a much stronger result, namely the noncommutative
analogue of Schreyer's Theorem [Sch] (cf. Theorem 3.7.13 in [AL] for
free modules over commutative polynomial algebras; Theorem 4.8 in
[Lev] for free modules over solvable polynomial
algebras):{\parindent=0pt \v5

{\bf 7.1. Theorem}  With respect to the left monomial   ordering  
$\prec_{s\hbox{-}\varepsilon}$ on $\B (\varepsilon )$ as defined 
above, the following statements hold.\par

(i) Let $s_{ij}$ be determined by $S_{\ell}(g_i,g_j)$, where $i<j$,
$\LM (g_i)=a^{\alpha (i)}e_s$ with $\alpha (i)=(\alpha_{i_1},\ldots
,\alpha_{i_n})$, and $\LM (g_j)=a^{\alpha (j)}e_s$ with $\alpha
(j)=(\alpha_{j_1},\ldots ,\alpha_{j_n})$. Then $\LM
(s_{ij})=a^{\gamma-\alpha (j)}\varepsilon_j$, where $\gamma
=(\gamma_1,\ldots ,\gamma_n)$ with each
$\gamma_k=\max\{\alpha_{i_k},\alpha_{j_k}\}$. \par

(ii) ${\cal S}$ is a left Gr\"obner basis of Syz$(\G )$, thereby
${\cal S}$ generates Syz$(\G )$.\par\QED} \v5

To go further, again let $\G=\{ g_1,\ldots ,g_t\}$ be  the left
Gr\"obner basis of $N$ produced by running {\bf Algorithm 1} with
the initial input data $ U=\{\xi_1,\ldots ,\xi_m\}$. Using the usual 
matrix notation for convenience, we have
$$\left (\begin{array}{l} \xi_1\\ \vdots\\ \xi_m\end{array}\right )=U_{m\times t}\left (\begin{array}{l} g_1
\\ \vdots\\ g_t\end{array}\right ) ,\quad
\left (\begin{array}{l} g_1\\ \vdots\\ g_t\end{array}\right
)=V_{t\times m}\left (\begin{array}{l} \xi_1
\\ \vdots\\ \xi_m\end{array}\right ) ,$$
where the $m\times t$ matrix $U_{m\times t}$ (with entries in $A$)
is obtained by the division by $\G$, and the $t\times m$ matrix
$V_{t\times m}$ (with entries in $A$) is obtained by keeping track
of the reductions during executing the WHILE loop of {\bf Algorithm
1}.  By Theorem 7.1, we may write Syz$(\G )=\sum^r_{i=1}A\mathcal
{S}_i$ with $\mathcal{S}_1,\ldots ,\mathcal{S}_r\in
L_1=\oplus_{i=1}^tA\varepsilon_i$;  and if
$\mathcal{S}_i=\sum^t_{j=1}f_{ij}\varepsilon_j$, then we  write
$\mathcal{S}_i$ as a $1\times t$ row matrix, i.e.,
$\mathcal{S}_i=(f_{i_1}~\ldots~ f_{it})$,  whenever matrix notation
is convenient in the according discussion. At this point, we note
also that all the $\mathcal{S}_i$ may be written down during
executing the WHILE loop of {\bf Algorithm 1} successively.
Furthermore, we write $D_{(1)},\ldots ,D_{(m)}$ for the rows of the
matrix $D_{m\times m}=U_{m\times t}V_{t\times m}-E_{m\times m}$
where $E_{m\times m}$ is the $m\times m$ identity matrix. The
following proposition is a noncommutative analogue of ([AL],
Theorem 3.7.6). {\parindent=0pt \v5

{\bf Proposition 7.2.}  With notation fixed above, the syzygy module
Syz$( U )$ of $ U =\{ \xi_1,\ldots ,\xi_m\}$ is generated by
$$\{ \mathcal{S}_1V_{t\times m},\ldots ,\mathcal{S}_{r}V_{t\times m},D_{(1)},\ldots ,D_{(m)}\} ,$$
where each $1\times m$ row matrix represents an element of the free
$A$-module $\oplus_{i=1}^mA\omega_i$. \vskip 6pt

{\bf Proof} Since $$0=\mathcal{S}_i\left (\begin{array}{l} g_1\\
\vdots\\ g_t\end{array}\right )=(f_{i_1}~\ldots ~f_{it})\left
(\begin{array}{l} g_1\\ \vdots\\ g_t\end{array}\right
)=(f_{i_1}~\ldots ~f_{it})V_{t\times m}\left (\begin{array}{l}
\xi_1\\ \vdots\\ \xi_m\end{array}\right ) ,$$ we have
$\mathcal{S}_iV_{t\times m}\in$ Syz$( U )$, $1\le i\le r$. Moreover, 
since
$$\begin{array}{rcl} D_{m\times m}\left (\begin{array}{l}
\xi_1\\ \vdots\\ \xi_m\end{array}\right )&=&(U_{m\times t}V_{t\times
m}-E_{m\times m})\left (\begin{array}{l} \xi_1\\ \vdots\\
\xi_m\end{array}\right )\\
\\
&=&U_{m\times t}V_{t\times m}\left (\begin{array}{l} \xi_1\\
\vdots\\ \xi_m\end{array}\right )-\left (\begin{array}{l} \xi_1\\
\vdots\\ \xi_m\end{array}\right )\\
\\
&=&U_{m\times t}\left (\begin{array}{l} g_1\\ \vdots\\
g_t\end{array}\right )-\left (\begin{array}{l} \xi_1\\ \vdots\\
\xi_m\end{array}\right )=\left (\begin{array}{l} \xi_1\\ \vdots\\
\xi_m\end{array}\right )-\left (\begin{array}{l} \xi_1\\ \vdots\\
\xi_m\end{array}\right )=0,\end{array}$$ we have $D_{(1)},\ldots
,D_{(r)}\in$ Syz$( U )$.}\par

On the other hand, if $H =(h_1~\ldots~h_m)$ represents the element
$\sum^m_{i=1}h_i\omega_i\in \oplus_{i=1}^mA\omega_i$ such that
$H\left (\begin{array}{l} \xi_1\\ \vdots\\ \xi_m\end{array}\right
)=0$, then $0=H U_{m\times t}\left (\begin{array}{l} g_1\\ \vdots\\
g_t\end{array}\right )$. This means $H U_{m\times t}\in$ Syz$(\G )$.
Hence, $H U_{m\times t}=\sum^r_{i=1}f_i\mathcal{S}_i$ with $f_i\in
A$, and it follows that $H U_{m\times t}V_{t\times
m}=\sum^r_{i=1}f_i\mathcal{S}_iV_{t\times m}$. Therefore,
$$\begin{array}{rcl} H&=&H +H U_{m\times t}V_{t\times m}-H U_{m\times t}V_{t\times m}\\
&=&H (E_m-U_{m\times t}V_{t\times
m})+\sum^r_{i=1}f_i\mathcal{S}_iV_{t\times m}\\
&=&-H D_{m\times m}+\sum^r_{i=1}f_i(\mathcal{S}_iV_{t\times
m}).\end{array}$$ This shows that every element of Syz$( U )$ is
generated by $\{ \mathcal{S}_1V_{t\times m},\ldots
,\mathcal{S}_{r}V_{t\times m},D_{(1)},\ldots ,D_{(m)}\} ,$ as
desired.\QED\v5

Next,  we recall the noncommutative version of Hilbert's syzygy
theorem for solvable polynomial algebras. For a constructive proof
of Hilbert's syzygy theorem by means of Gr\"obner bases respectively
in  the commutative case and the noncommutative case, we refer to
(Corollary 15.11 in  [Eis]) and (Section 4.4 in [Lev]).
{\parindent=0pt\v5

{\bf 7.3. Theorem} Let $A=K[a_1,\ldots ,a_n]$ be an arbitrary
solvable polynomial algebra with  admissible system $(\B ,\prec )$.
Then every finitely generated (left) $A$-module $M$ has a free
resolution
$$0~\mapright{}{}~L_s~\mapright{}{}~L_{s-1}~\mapright{}{}\cdots ~
\mapright{}{}~L_0~\mapright{}{}~M~\mapright{}{}~ 0$$ where each
$L_i$ is a free $A$-module of finite rank and $s\le n$. It follows
that $M$ has projective dimension p.dim$_AM\le s$, and that $A$ has
global homological dimension gl.dim$A\le n$. \par\QED}\v5

Now, we are ready to reach the goal of this section.
{\parindent=0pt\v5

{\bf 7.4. Theorem} Let $A=K[a_1,\ldots ,a_n]$ be the solvable
polynomial algebra fixed in the beginning of this section, and let
$L_0=\oplus_{i=1}^mAe_i$ be a filtered free $A$-module with the
filtration $FL_0=\{ F_qL_0\}_{q\in\NZ}$ such that $d_{\rm
fil}(e_i)=b_i$, $1\le i\le m$. If $N=\sum^s_{i=1}A\xi_i$ is a
finitely generated submodule of $L_0$ and the quotient module
$M=L_0/N$ is equipped with the filtration
$FM=\{F_qM=(F_qL_0+N)/N\}_{q\in\NZ}$. Then $M$ has a minimal
filtered free resolution of length $d\le n$ (in the sense of
Definition 6.2):
$${\cal L}_{\bullet}\quad\quad 0~\mapright{}{}~L_d~\mapright{\varphi_{q}}{}~
\cdots~\mapright{\varphi_2}{}~L_1~\mapright{\varphi_1}{}~L_0~\mapright{\varphi_0}{}~M~\mapright{}{}~0$$
which can be constructed by implementing the following
procedures:}\par

{\bf Procedure 1.} Fix a graded left monomial ordering
$\prec_{e\hbox{-}gr}$ on the $K$-basis $\BE$ of $L_0$ (see Section 
3), and run {\bf Algorithm 1} with the initial input data $ U =\{ 
\xi_1,\ldots ,\xi_s\}$ to compute a left Gr\"obner basis $\G=\{ 
g_1,\ldots ,g_z\}$ for $N$, so that $N$ has the standard basis $\G$ 
with respect to the induced filtration $FN=\{ F_qN=N\cap 
F_qL_0\}_{q\in\NZ}$ (Theorem 4.8).\par

{\bf Procedure 2.} Run {\bf Algorithm 2} (presented in Section 5)
with the initial input data $E=\{ e_1,\ldots ,e_m\}$ and $\G=\{
g_1,\ldots ,g_z\}$ to compute a subset $\mathscr{E}_0'=\{
e_{i_1},\ldots ,e_{i_{m'}}\}\subset \mathscr{E}_0=\{ e_1,\ldots
,e_m\}$ and a subset $V=\{ v_1,\ldots ,v_t\}\subset N\cap L_0'$ such
that there is a strict filtered isomorphism $L_0'/N'=M'\cong M$,
where $L_0'=\oplus_{q=1}^{m'}Ae_{i_q}$ with $m'\le m$ and
$N'=\sum^{t}_{k=1}Av_k$, and such that $\{\OV e_{i_1},\ldots ,\OV
e_{i_{m'}}\}$ is a minimal F-basis of $M$ with respect to the
filtration $FM$.\par

For convenience, after accomplishing Procedure 2 we may assume that
$\mathscr{E}_0=\mathscr{E}_0'$, $ U =V$ and $N=N'$. Accordingly we 
have the short exact sequence $0~\mapright{}{}~ N~\mapright{}{}~ 
L_0~\mapright{\varphi_0}{}~M~\mapright{}{}  ~0$ such that $\varphi_0 
(\mathscr{E}_0)=\{ \OV e_1,\ldots ,\OV e_m\}$  is a minimal  F-basis 
of $M$ with respect to the filtration $FM$. \par

{\bf Procedure 3.} With the initial input data $ U =V$, implements 
the procedures presented in Theorem 5.4 to compute a minimal 
standard basis $W=\{ \xi_{j_1},\ldots ,\xi_{j_{m_1}}\}$ for $N$ with 
respect to the induced filtration $FN$.
\par

{\bf Procedure 4.} Computes a generating set $U_1  =\{\eta_1,\ldots 
,\eta_{s_1}\}$ of $N_1=$ Syz$(W)$ in the free $A$-module 
$L_1=\oplus_{i=1}^{m_1}A\varepsilon_i$ by running {\bf Algorithm 1} 
with the initial input data $W$ and using Proposition 7.2.\par

{\bf Procedure 5.} Construct the strict filtered exact sequence
$$0~\mapright{}{}~ N_1~\mapright{}{}~L_1~\mapright{\varphi_1}{}~
L_0~\mapright{\varphi_0}{}~M~\mapright{}{}  ~0$$ where the
filtration $FL_1$ of $L_1$ is constructed by setting $d_{\rm
fil}(\varepsilon_k)=d_{\rm fil}(\xi_{j_k})$, $1\le k\le m_1$, and
$\varphi_1$ is defined by setting
$\varphi_1(\varepsilon_k)=\xi_{j_k}$, $1\le k\le m_1$. If $N_1\ne
0$, then, with the initial input data $ U =U_1$,  repeat Procedure 3 
-- Procedure 5 for $N_1$ and so on.\par

By Theorem 6.3, a minimal filtered free resolution ${\cal
L}_{\bullet}$ of $M$ gives rise to a minimal graded free resolution
$G({\cal L}_{\bullet})$ of $G(M)$. Since $G(A)=K[\sigma (a_1),\ldots
,\sigma (a_n)]$ is a solvable polynomial algebra by Theorem 2.4, it
follows from  Theorem 7.3 that $G({\cal L}_{\bullet})$ terminates at
a certain step, i.e., Ker$G(\varphi_d)=0$ for some $d$. But
Ker$G(\varphi_d)=G(\hbox{Ker}\varphi_d)$ by Proposition 5.1, where
Ker$\varphi_d$ has the filtration induced by $FL_d$, thereby
$G(\hbox{Ker}\varphi_d)=0$. Consequently Ker$\varphi_d=0$ by
classical filtered modules theory, thereby a minimal finite filtered
free resolution of length $d\le n$ is achieved for $M$. \v5

\v5 \centerline{References}
\parindent=1.6truecm

\re{[AL1]} J. Apel and W. Lassner, An extension of Buchberger's
algorithm and calculations in enveloping fields of Lie algebras.
{\it J. Symbolic Comput}., 6(1988), 361--370.

\item{[AL2]} W. W. Adams and P. Loustaunau, {\it An Introduction to Gr\"obner Bases}.
Graduate Studies in Mathematics, Vol. 3. American Mathematical
Society, 1994.

\item{[AVV]} M. J. Asensio, M. Van den Bergh and F. Van Oystaeyen,
A new algebraic approach to microlocalization of filtered rings.
Trans. Amer. math. Soc, 2(316)(1989), 537-555.

\item{[Bu1]} B. Buchberger, {\it Ein Algorithmus zum Auffinden der
Basiselemente des Restklassenringes nach einem nulldimensionalen
polynomideal}. PhD thesis, University of Innsbruck, 1965.

\item{[Bu2]} B. Buchberger, Gr\"obner bases: ~An algorithmic method
in polynomial ideal theory. In: {\it Multidimensional Systems
Theory} (Bose, N.K., ed.), Reidel Dordrecht, 1985, 184--232.

\item{[BW]} T.~Becker and V.~Weispfenning, {\it Gr\"obner Bases}.
Springer-Verlag, 1993.

\item{[CDNR]} A. Capani, G. De Dominicis, G. Niesi, and L. Robbiano,
Computing minimal finite free resolutions. {\it Journal of Pure and
Applied Algebra}, (117\& 118)(1997), 105 -- 117.

\item{[DGPS]} W. Decker, G.-M. Greuel, G. Pfister, and H. Sch{\"o}nemann, \newblock {\sc Singular} {3-1-3}
--- {A} computer algebra system for polynomial computations.
\newblock {http://www.singular.uni-kl.de}(2011).

\item{[Eis]} D. Eisenbud, {\it Commutative Algebra with a View
Toward to Algebraic Geometry}, GTM 150. Springer, New York, 1995.

\item{[Fr\"ob]} R. Fr\"oberg, {\it An Introduction to Gr\"obner Bases}. Wiley, 1997.

\item{[Gal]} A. Galligo, Some algorithmic questions on ideals of
differential operators. {\it Proc. EUROCAL'85}, LNCS 204, 1985,
413--421.

\item{[Gol]} E. S. Golod, Standard bases and homology, in: {\it Some
Current Trends in Algebra}. (Varna, 1986), Lecture Notes in
Mathematics, 1352, Springer-Verlag, 1988, 88-95.

\item{[Kr]} H. Kredel, {\it Solvable Polynomial Rings}. Shaker-Verlag, 1993.

\item{[KR]} M. Kreuzer, L. Robbiano, {\it Computational Commutative Algebra 2}. Springer, 2005.

\item{[K-RW]} A.~Kandri-Rody and V.~Weispfenning, Non-commutative
Gr\"obner bases in algebras of solvable type. {\it J. Symbolic
Comput.}, 9(1990), 1--26.

\item{[Lev]} V. Levandovskyy, {\it Non-commutative Computer Algebra for
Polynomial Algebra}: {\it Gr\"obner Bases, Applications and
Implementation}. Ph.D. Thesis, TU Kaiserslautern, 2005.

\item{[Li1]} H. Li, {\it Noncommutative Gr\"obner Bases and
Filtered-Graded Transfer}. Lecture Notes in Mathematics, Vol. 1795,
Springer-Verlag, 2002.

\item{[Li2]} H. Li, A Constructive characterization of solvable polynomial algebras.
arXiv:1212.5988 [math.RA].

\item{[Li3]} H. Li, Computation of minimal  graded free resolutions over $\NZ$-graded solvable polynomial algebras.
arXiv:1401.???? [math.RA].

\item{[Li4]} H. Li, {\it Gr\"obner Bases in Ring Theory}. World Scientific Publishing Co., 2011.

\item{[LS]} H. Li and C. Su, On (de)homogenized Gr\"obner bases. {\it
Journal of Algebra, Number Theory}: {\it Advances and Applications},
3(1)(2010), 35--70

\item{[LVO]} H. Li and F. Van Oystaeyen, {\it Zariskian
Filtrations}. $K$-Monograph in Mathematics, Vol.2. Kluwer Academic
Publishers, 1996.

\item{[LW]} H. Li and  Y. Wu, ~Filtered-graded transfer of
Gr\"obner basis computation in solvable polynomial algebras. {\it
Communications in Algebra}, 1(28)(2000), 15--32.

\item{[Sch]} F.O. Schreyer, {\it Die Berechnung von Syzygien mit
dem verallgemeinerten Weierstrasschen Divisionsatz}. Diplomarbeit,
Hamburg, 1980.

\end{document}